\documentclass[final]{siamltex}
\usepackage[dvips]{graphics}
\usepackage{amsmath}
\input{amssym.def}
\input{amssym.tex}
\title{Modeling of a curvilinear planar crack with a curvature-dependent surface tension\thanks{This work was
supported by Award No. KUS-C1-016-04, made by King Abdullah University of Science and Technology (KAUST)}}

\author{A. Y. ZEMLYANOVA, J.R. WALTON\thanks{Department of Mathematics, IAMCS, Texas A\&M University. 
(azem@math.tamu.edu). Questions, comments, or corrections
to this document may be directed to that email address.}}

\def\theequation{\thesection.\arabic{equation}}

\date{}

\newcommand{\I}{\mathop{\rm Im}\nolimits}
\newcommand{\R}{\mathop{\rm Re}\nolimits}

\newcommand{\beqa}{\begin{eqnarray}}
\newcommand{\eeqa}[1]{\label{#1}\end{eqnarray}}
\newcommand{\bequ}{\begin{equation}}
\newcommand{\eequ}[1]{\label{#1}\end{equation}}

\newcommand{\beq}{\begin{equation}}
\newcommand{\eeq}{\end{equation}}
\newcommand{\barr}{\begin{eqnarray}}
\newcommand{\earr}{\end{eqnarray}}
\newcommand{\beqn}{\begin{equation*}}
\newcommand{\eeqn}{\end{equation*}}
\newcommand{\barrn}{\begin{eqnarray*}}
\newcommand{\earrn}{\end{eqnarray*}}

\newcommand{\bfa}[1]{\mathbf{#1}}     %
\newcommand{\divv}{\mbox{div}}        %
\newcommand{\tgamma}{{\tilde{\gamma}}}%
\newcommand{\jump}[1]{[\![ #1 ]\!]}

\begin{document}
\maketitle

\noindent

\begin{abstract}
An approach to modeling fracture incorporating interfacial mechanics is applied to the example of a curvilinear plane strain crack. The classical Neumann boundary condition is augmented with curvature-dependent surface tension. It is shown that the considered model eliminates the integrable  crack-tip stress and strain singularities of order $1/2$ present in the classical linear fracture mechanics solutions, and also leads to the sharp crack opening that is consistent with empirical observations. Unlike for the case of a straight crack, for a general curvilinear crack some components of the stresses and the derivatives of the displacements may still possess weaker singularities of a logarithmic type. Generalizations of the present study that lead to complete removal of all crack-tip singularities, including logarithmic, are the subject of a future paper. 

\end{abstract}



\pagestyle{myheadings}
\thispagestyle{plain}
\markboth{A. Y. ZEMLYANOVA, J.R. WALTON}{Modeling of a curvilinear planar crack...}

\setcounter{equation}{0}

\section{Introduction}

A better understanding of fracture processes is of utmost importance for applications. The standard approach to study the behavior and propagation of fracture is based within the framework of linear elastic fracture mechanics (LEFM). This approach has been successfully applied to countless problems and has a vast literature. However, LEFM contains a well-known internal inconsistency. While it is based on the assumption that the stresses and the strains remain small everywhere in the body, LEFM predicts that the stresses and the strains possess an integrable power singularity of the order $1/2$ at the crack tips. Moreover, in the case of an interface crack between two dissimilar materials there is an additional oscillating singularity which results in the non-physical interpenetration of the crack faces near the crack tips.

Several models have been proposed to eliminate this internal inconsistency. One of the most common approaches, studied by many authors, is to introduce two-dimensional cohesive zones or three-dimensional process zones near the crack tips. Despite the obvious advantages of this approach it has shortcomings of both theoretical and practical nature, such as difficulty in specifying physically valid constitutive response functions.

Since the fracture occurs within nano- or molecular scale processes, it has been argued that it is impossible to describe the fracture effectively within the context of continuum mechanics. Thus, there is a growing literature dedicated to the modeling of fracture using atomistic and lattice based approaches \cite{Abraham2001, AbrahamEtal1977, AbrahamGao2000, FinebergGross1991, HollandMarder1997, MarderGross1995, SlepyanEtAl1999, SwadenerEtAl2002}. The accuracy of these methods largely depends on the precise description of intermolecular forces which is difficult to do for liquids and solids \cite{HollandMarder1997}. This approach may also present some computational challenges. Several atom-to-continuum models have been proposed as well. One of the most extensively studied methods of this type, in the context of finite element method (FEM) approximations to continuum
models, is the quasi-continuum method introduced by Tadmor et al in 1996 \cite{TadmoretAl1996}. Based on an atomistic view of material behavior, its continuum aspect comes from the fact that the FEM is based on energy minimization. A different type of atom-to-continuum modeling is a recently proposed approach by Xiao and Belytschko in \cite{XiaoBelytschko2004} which involves the introduction of bridging domains between regions modeled using bulk (continuum) descriptions of material behavior and regions modeled using atomistic descriptions of material behavior. Both of these approaches involve adjustable parameters, such as choosing the size and location of domains over which  potentials acting at different length scales are blended, that one needs to fit for every particular application.

In a series of works \cite{FometheMaugin1998, GurtinPodio1996, GurtinPodio1998, GurtinShvartsman1997, MauginTrimarco1995, Spector2002} the modeling of fracture is based on the introduction of the concept of configurational forces. The configurational forces were first introduced as derivatives of energy by Eshelby \cite{Eshelby1956} and were further studied by Gurtin \cite{Gurtin1995}. The configurational forces are described in \cite{Gurtin1995} as the intrinsic properties of a body's material micro-structure which have their own force balance.

Another continuum approach, peridynamics \cite{HaBobaru2010a, HaBobaru2011, Silling2000}, addresses the modeling of the spontaneous formation and dynamic propagation of a fracture, for which the classical continuum mechanics approach is ill-suited as well. In peridynamics, material points are assumed to be connected via elastic (linear or nonlinear) bonds which have a critical relative elongation. This model does not distinguish between the points in a body where a discontinuity in displacement or its derivatives may be located. The essence of the mathematical aspect of this nonlocal model is that integration rather than differentiation is used to compute the force on a material particle.

One of the most convenient and successful modifications to LEFM incorporates surface excess properties on the crack boundary into the modeling the fracture. It is well-known that the molecules on the boundary between two different phases are subjected to different forces than molecules in the bulk of the material. This model accounts for this fact by introducing a dividing surface endowed with surface excess properties \cite{SlatteryEtAl2004}. This theory has been first applied to fracture mechanics in \cite{OhWaltonSlattery2006} and developed in more details on the example of a straight mode-I crack in \cite{SendovaWalton2010}. The advantages of this model are that, first, it relies on the physically valid assumption that the forces which act on the boundary of the object are different from the forces in the bulk and a curvature-dependent surface tension is present along the interface between two phases. It has been observed also that the behavior of materials is in accordance with linear elasticity everywhere except very small neighborhoods around the crack tips which provides additional justification for using the equations of linear elasticity to describe the bulk of material. Hence, since the standard equations of linear elasticity are assumed for the bulk of the material and only the boundary conditions are changed, it follows that a wide variety of mathematical approaches and techniques used in linear elasticity are still applicable to this model. Finally, this model predicts bounded stresses and strains in the case of a straight mode-I crack, thus removing the inconsistency present in LEFM. This model also bypasses the introduction of artificial conditions on the boundary and avoids ad hoc choices of parameters that are present in most other models. 

While the body of research discussed above is dedicated to modification of the LEFM in order to remove the power singularity of the order $1/2$ at the crack tips, it is necessary to note that most of these attempts, with rare exception, consider straight cracks only. For practical applications, however, it is important to study curvilinear cracks. In this paper we extend the curvature dependent surface tension model \cite{SendovaWalton2010} to the curvilinear smooth plane strain crack of essentially arbitrary shape. The most significant conclusion of the paper is that the analysis of a curvilinear (non-zero curvature) crack is significantly different from the case of a straight crack. Even though just as in the case of the straight crack, the stresses and strains do not have an integrable power singularity of the order $1/2$, thus removing the main inconsistency of LEFM, the shear stress and one of the derivatives of displacements may possess a weaker logarithmic singularity while the tensile stress and the derivative of the other displacement remain bounded. Thus, the curvature of the crack plays an important role in the singular behavior of the stresses and strains. This suggests that controlling mean curvature alone may not be enough to force the stresses and strains to be bounded. We conjecture that the present theory must be extended through the addition of higher surface gradient effects in the model for stress-deformation behavior within the fracture surfaces in order to remove the logarithmic singularities predicted by the curvature dependent surface tension model considered here. Such considerations are beyond the scope of the present contribution and are the subject of an ongoing research effort.

The complex analysis approach is used for the solution of the stated problem. Functions of complex variable have been first applied to the solution of two-dimensional problems of elasticity by Kolosov in 1909. The methods of complex analysis have been further developed by Muskhelishvili \cite{Mus1963} who inspired a large school of research in this area. The complex analysis approach is an extremely powerful technique and can be applied to a large number of problems of two-dimensional elasticity. The main idea consists in expressing the stresses and the displacements in elastic body through two analytic functions, thus, reducing the problem of elasticity from the solution of partial differential equations to the solution of boundary problem for these two analytic functions. This boundary problem is often reduced further either to Riemann-Hilbert type of problem or to integral equations of some kind. Some excellent applications of methods of complex analysis in elasticity can be found in \cite{England1971}, \cite{Mus1963}, \cite{Sokolnikoff1956} and others. 

The solution of the problem proceeds in the following way. First, the boundary condition of the type \cite{SendovaWalton2010} is stated. The curvature of the crack dictates a slight modification to the form of the surface tension considered in \cite{SendovaWalton2010}. In particular, it is reasonable to assume zero surface tension in the unloaded configuration of the body. Next, this condition is linearized and written in terms of the tensile and the shear stresses and derivatives of the displacements. Muskhelishvili's formulas \cite{Mus1963} are applied then to reduce this condition to the boundary value problem for two analytic functions (complex potentials) and their derivatives. On the next step, the boundary value problem is reduced to the system of two Cauchy-type singular integro-differential equations with the help of Savruk's integral representations of complex potentials \cite{Savruk1981}. This system is further reduced to the system of two weakly singular integral equations in the spirit of \cite{MikhPros1986}. It is necessary to note that the mathematical techniques used to solve the problem are significantly different from those in \cite{SendovaWalton2010}. Finally, the numerical procedure is developed to solve the system of singular integro-differential equations and the numerical results are presented.

\setcounter{equation}{0}

\section{Model with a curvature dependent surface tension}

Consider the differential and jump momentum conditions in the deformed configuration in the absence of inertial and gravitational effects \cite{Slattery1990}:
\begin{align}
\mbox{div}(\bfa{T}) &=0,
\label{2_1}\\
\grad_{(\zeta)}\tgamma+2\tgamma H\bfa{n} +\jump{\bfa{T}}\bfa{n} &=\bfa{0},
\label{2_2}
\end{align}
where $\bfa{T}$ is the Cauchy stress tensor, $\bfa{n}$ is the unit normal to the fracture surface $\zeta$ pointing into the bulk material, $H=-\frac12\divv_{(\zeta)}{\bfa{n}}$ is the mean curvature, subscript $\ldots_{(\zeta)}$ is used to emphasize that the quantity is associated with a surface of the crack, $\grad_{(\zeta)}$ and $\divv_{(\zeta)}$ denote the surface gradient and the surface divergence correspondingly with respect to the given parametrization of the surface and can be derived using the formulas in \cite{Slattery1990}, and the double brackets $\jump{\ldots}$ denote the jump of the quantity enclosed across the boundary of the crack.

Assume that $\tgamma$, which will be called surface tension, depends linearly on the difference between the curvature of the deformed crack and the curvature of the crack in the unloaded configuration:

\begin{equation}
\tilde{\gamma}=\gamma_1(\divv_{(\zeta)}{\bfa{n}}-\divv_{(\zeta_0)}{\bfa{n}_0}),
\label{2_3}
\end{equation}
where the subindex ``0" denotes the parameters in the initial configuration of the crack.

\begin{figure}[ht]
	\centering
		\scalebox{0.5}{\includegraphics{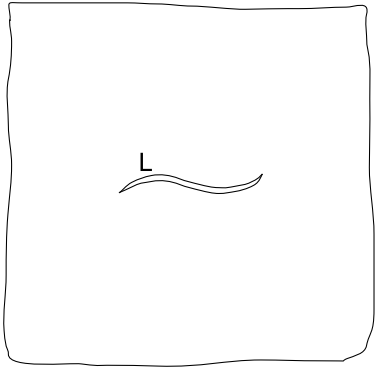}}
			\caption{An infinite plate $S$ with a crack $L$.}
	\label{fig0}
\end{figure}

Since the crack $L$  in an undeformed configuration is a Jordan arc (fig. \ref{fig0}), it is natural to take the parametric equations of the crack $L$ to be $t(s)=x_1(s)+ix_2(s)$, where the parameter $s$ is an arc length, $s\in [0,l]$. Assume that the function $t(s)$ has continuous derivatives up to the fourth order. This assumption is made for simplicity and can be somewhat relaxed. Since the boundary of the crack is free from the external stresses and, hence, the jump of the stresses on each of the boundaries of the crack is simply equal to the stresses in the cracked plate, it is possible to rewrite the condition (\ref{2_2}) in terms of the tensile and shear stress in the following form:
\begin{equation}
\sigma_n+i\tau_n=i(\grad_{(\zeta)}\tgamma-\tgamma\bfa{n}\divv_{(\zeta)}\bfa{n})\overline{t'(s)},
\label{2_4}
\end{equation}
where $\sigma_n$ and $\tau_n$ are the tensile and shear stresses on the crack boundary.
Furthermore, it is easy to see that in the deformed configuration of the crack the quantity $\divv_{(\zeta)}\bfa{n}$ can be expressed as 
\begin{equation}
\divv_{(\zeta)}\bfa{n}=\frac{\ddot{X}_2\dot{X}_1-\ddot{X}_1\dot{X}_2}{(\dot{X}_1^2+\dot{X}_2^2)^{3/2}},
\label{2_5}
\end{equation}
where $X_1(s)+iX_2(s)=(x_1(s)+u_1(s))+i(x_2(s)+u_2(s))$ are the equations of the crack in the deformed state and $(u_1+iu_2)(s)$ are the displacements. Linearizing the condition (\ref{2_5}) under the assumption that the derivatives of the displacements are small, one obtains
\begin{align}
\divv_{(\zeta)}\bfa{n} &=\kappa_0(s)+(\ddot{x}_2(s)-3\dot{x}_1(s)\kappa_0(s))\frac{du_1}{ds}-(\ddot{x}_1(s)+3\dot{x}_2(s)\kappa_0(s))\frac{du_2}{ds}
\notag\\
&-\dot{x}_2(s)\frac{d^2u_1}{ds^2}+\dot{x}_1(s)\frac{d^2u_2}{ds^2},
\label{2_6}
\end{align}
where $\kappa_0(s)$ is the curvature of the crack in the undeformed configuration.
Substituting (\ref{2_6}) into (\ref{2_4}) leads to the following boundary condition
\begin{align}
(\sigma_n+i\tau_n)^{\pm}(s)&=m_1(s)\frac{d}{ds}(u_1+iu_2)^{\pm}(s)+m_2(s)\frac{d}{ds}{(u_1-iu_2)}^{\pm}(s)+m_3(s)\frac{d^2}{ds^2}(u_1+iu_2)^{\pm}(s)
\notag\\
&+m_4(s)\frac{d^2}{ds^2}{(u_1-iu_2)}^{\pm}(s)+i\gamma_1\I \frac{d^3}{ds^3}\left(\overline{t'(s)}(u_1+iu_2)^{\pm}(s)\right),\,\,\,\,s\in [0,l],
\label{2_7}\\
m_1(s)&=-0.5\gamma_1(\overline{t'''(s)}+2i\overline{t''(s)}\kappa_0(s)+3i\overline{t'(s)}\kappa'_0(s)+3\overline{t'(s)}\kappa_0^2(s)),
\notag\\
m_2(s)&=0.5\gamma_1({t'''(s)}-4i{t''(s)}\kappa_0(s)-3i{t'(s)}\kappa'_0(s)-3{t'(s)}\kappa^2_0(s)),
\notag\\
m_3(s)&=-2i{\gamma_1}\overline{t'(s)}\kappa_0(s),\,\,\,\,m_4(s)=-i{\gamma_1}{t'(s)}\kappa_0(s).
\notag
\end{align}
where ``$+$" corresponds to the stresses and the displacements on the left-hand side of the crack with respect to the direction chosen, and ``$-$" corresponds to the same functions on the right-hand side of the crack.
Thus, the equation (\ref{2_7}) is valid on the both surfaces of the crack.

\section{Complex potentials and their integral representations}
It is well known \cite{Mus1963} that the tensile $\sigma_n$ and the shear $\tau_n$ components of the stress vector acting on the tangent line to the curve  $L$ from the side of the positive normal and the derivative of the displacement vector $u_1+iu_2$ at the points of the curve $L$ can be obtained from the Muskhelishvili's complex potentials $\Phi(z)$, $\Psi(z)$ by the following formulas \cite{Savruk1981}:\\
\begin{align}
(\sigma_n+i\tau_n)(t)&=\Phi(t)+\overline{\Phi(t)}+\frac{\overline{dt}}{dt}(t\overline{\Phi'(t)}+\overline{\Psi(t)}),
\label{3_1}\\
2\mu\frac{d}{dt}(u_1+iu_2)(t)&=\kappa\Phi(t)-\overline{\Phi(t)}-\frac{\overline{dt}}{dt}(t\overline{\Phi'(t)}+\overline{\Psi(t)}), \,\, t\in L.
\label{3_2}
\end{align}
Here $\Phi(z)$ and $\Psi(z)$ are analytic functions outside of the crack, $\mu$ is the shear modulus of the material and $\kappa=(3-\nu)/(1+\nu)$ for the case of plane stress and $\kappa=3-4\nu$ for plane strain, where $\nu$ is the Poisson's ratio.

Consider the infinite plate $S$ with a crack $L$ under the action of the given at infinity stresses. Assume that the stresses $\sigma_n+i\tau_n$ and the derivatives of the displacements $\frac{d}{dt}(u_1+iu_2)$ in the plate have jump discontinuities on the curve $L$. Then the complex potentials $\Phi(z)$ and $\Psi(z)$ can be taken in the form \cite{Savruk1981}:\\
\begin{align}
\Phi(z)&=\Gamma+\frac{(\kappa+1)^{-1}}{2\pi}\int_L\frac{g'(t)dt}{t-z}+\frac{(\kappa+1)^{-1}}{\pi i}\int_{L}\frac{q(t)dt}{t-z},
\label{3_3}\\
\Psi(z)&=\Gamma'+\frac{(\kappa+1)^{-1}}{2\pi}\int_L\left( \frac{\overline{g'(t)dt}}{t-z}-\frac{\bar{t}g'(t)dt}{(t-z)^2}\right)+\frac{(\kappa+1)^{-1}}{\pi i}\int_{L}\left(\frac{\kappa\overline{q(t)dt}}{t-z}-\frac{\bar{t}q(t)dt}{(t-z)^2}\right),
\notag\\
z&\in S,\,\,\,\Gamma=(\sigma_1^{\infty}+\sigma_2^{\infty})/4, \,\, \Gamma'=(\sigma_2^{\infty}-\sigma_1^{\infty})e^{-2i\alpha}/2,
\notag
\end{align}
where $\sigma_1^{\infty}$ and $\sigma_2^{\infty}$ are the principal tensile stresses at infinity acting in the directions which constitute the angles $\alpha$ and $\alpha+\pi/2$ with the positive direction of $x$-axis correspondingly, and $g'(t)$ and $q(t)$ are the unknown functions which can be expressed through the jump of the displacements $u_1+iu_2$ and the jump of the stresses $\sigma_n+i\tau_n$ on the line $L$ by the formulas:\\
\begin{align}
g'(t)&=-2i\mu\frac{d}{dt}\left((u_1+iu_2)^+(t)-(u_1+iu_2)^-(t)\right), \,\, t\in L,
\label{3_4}\\
q(t)&=\left((\sigma_n+i\tau_n)^+(t)-(\sigma_n+i\tau_n)^-(t)\right)/2, \,\, t\in L.
\notag
\end{align}

Hence, the stressed state of the cracked plate is described by the complex potentials (\ref{3_3}) which contain two unknown functions $g'(t)$, $q(t)$ defined on the crack $L$. We will look for these functions in the class of functions satisfying the H\"older condition on the curve $L$. This choice guarantees \cite{Gakhov1990} the existence of all principal and limit values of the integrals of the Cauchy type in the formulas (\ref{3_3}). Also assume that the formula (\ref{3_2}) can be differentiated by $s$ twice. This will be proved rigorously later.

Using (\ref{3_1}), (\ref{3_2}) and (\ref{3_3}) it is possible to express the stresses $(\sigma_n+i\tau_n)$ and the derivatives of the displacements $\frac{d}{ds_0}(u_1+iu_2)$ in the form:
$$
(\sigma_n+i\tau_n)^{\pm}(s_0)=\pm q(s)+\frac{(\kappa+1)^{-1}}{2\pi}\left[\int_0^l\left(\frac{2}{s-s_0}+k_1(s,s_0) \right)g'(s)ds+\int_0^lk_2(s,s_0)\overline{g'(s)}ds\right.
$$
\begin{equation}
\label{3_5}
\left.-2i\int_0^l\left(-\frac{\kappa-1}{s-s_0}+k_3(s,s_0) \right)q(s)ds+2i\int_0^l k_2(s,s_0)\overline{q(s)}ds\right]
+2\R \Gamma +\overline{\Gamma'(t'(s_0))^2};
\end{equation}
\begin{equation}
\label{3_6}
2\mu\frac{d}{ds_0}(u_1+iu_2)^{\pm}(s_0)=t'(s_0)\Omega^{\pm}(s_0)+(\kappa\Gamma-\overline{\Gamma})t'(s_0)-\overline{\Gamma't'(s_0)},
\end{equation}
where the functions $\Omega^{\pm}(s_0)$ have the following form
\begin{align}
\Omega^{\pm}(s_0)&=\pm\frac{i(\kappa+1)}{2}g'(s_0)+\frac{(\kappa+1)^{-1}}{2\pi}\left[\int_0^l\left(\frac{\kappa-1}{s-s_0}+k_4(s,s_0) \right)g'(s)ds\right.
\label{3_7}\\
&-\int_0^lk_2(s,s_0)\left.\overline{g'(s)}ds-2i\int_0^l\left(\frac{2\kappa}{s-s_0}+\kappa k_1(s,s_0) \right)q(s)ds
-2i\int_0^l k_2(s,s_0)\overline{q(s)}ds\right],
\notag
\end{align}
where $g'(s)=g'(t(s))$, $q(s)=q(t(s))$ and $k_j(s,s_0)$, $j=1,\,2,\,3,\,4$, denote the regular kernels:
\begin{align}
k_1(s,s_0)&=-\frac{2}{s-s_0}+\frac{t'(s)}{t(s)-t(s_0)}+\frac{t'(s)}{\overline{t(s)}-\overline{t(s_0)}}\frac{\overline{t'(s_0)}}{t'(s_0)};
\notag\\
k_2(s,s_0)&=\frac{\overline{t'(s)}}{\overline{t(s)}-\overline{t(s_0)}}-\frac{(t(s)-t(s_0))\overline{t'(s)}}{(\overline{t(s)}-\overline{t(s_0)})^2}\frac{\overline{t'(s_0)}}{t'(s_0)};
\label{3_8}\\
k_3(s,s_0)&=\frac{\kappa-1}{s-s_0}+\frac{t'(s)}{t(s)-t(s_0)}-\frac{\kappa t'(s)}{\overline{t(s)}-\overline{t(s_0)}}\frac{\overline{t'(s_0)}}{t'(s_0)};
\notag\\
k_4(s,s_0)&=-\frac{\kappa-1}{s-s_0}+\frac{\kappa t'(s)}{t(s)-t(s_0)}-\frac{t'(s)}{\overline{t(s)}-\overline{t(s_0)}}\frac{\overline{t'(s_0)}}{t'(s_0)}.
\notag
\end{align}

It is possible to rewrite the boundary conditions (\ref{2_7}) in the following form:
\begin{align}
(\sigma_n+i\tau_n)^{\pm}(s_0)&=-\frac{\gamma_1}{2\mu}\kappa_0(s_0)\left[\kappa_0(s_0)\R \Omega^{\pm}(s_0)-\I\frac{d}{ds_0}\Omega^{\pm}(s_0)\right]
\notag\\
&-i\frac{\gamma_1}{2\mu}\frac{d}{ds_0}\left[\kappa_0(s_0)\R \Omega^{\pm}(s_0)-\I\frac{d}{ds_0}\Omega^{\pm}(s_0)\right] +f(s_0),
\label{3_9}
\end{align}
where the functions $\Omega^{\pm}(s_0)$ have been introduced above, and
\begin{align*}
f(s_0)&=\gamma_1(4\mu)^{-1}\left[m_1(s_0)((\kappa\Gamma-\overline{\Gamma})t'(s_0)-\overline{\Gamma't'(s_0)})+m_2(s_0)((\kappa\overline{\Gamma}-\Gamma)\overline{t'(s_0)}-\Gamma't'(s_0)) \right.\\
&+m_3(s_0)((\kappa\Gamma-\overline{\Gamma})t''(s_0)-\overline{\Gamma't''(s_0)})+m_4(s_0)((\kappa\overline{\Gamma}-\Gamma)\overline{t''(s_0)}-\Gamma't''(s_0))\\
&+\left.2i\I\{\overline{t'(s_0)}((\kappa\Gamma-\overline{\Gamma})t'''(s_0)-\overline{\Gamma't'''(s_0)})\}\right]
-2\R\Gamma-\overline{\Gamma}'(\overline{t'(s_0)})^2.
\end{align*}

The boundary conditions (\ref{3_9}) must be combined with an additional condition providing that the displacements are single-valued on the crack $L$ which can be expressed in the simple form as:
\begin{equation}
\label{3_9a}
\int_0^lg'(s)t'(s)ds=0.
\end{equation}

Observe that since $(\sigma_n+i\tau_n)^{\pm}(s_0)$ and $\Omega^{\pm}(s_0)$ are allowed to have at most integrable power singularities at the crack tips, it follows immediately that $\kappa_0(s_0)\R \Omega^{\pm}(s_0)-\I\frac{d}{ds_0}\Omega^{\pm}(s_0)$ must be bounded at the the crack tips, and hence, the stress $\sigma_n^{\pm}$ and the derivative of the displacement $du_2^{\pm}/ds$ are bounded. Additionally, boundedness of $\kappa_0(s_0)\R \Omega^{\pm}(s_0)-\I\frac{d}{ds_0}\Omega^{\pm}(s_0)$ makes it possible \cite{MikhPros1986} to differentiate the Cauchy singular integrals in (\ref{3_9}). Substitute (\ref{3_5}), (\ref{3_7}) into (\ref{3_9}) to obtain the following system of equations:
\begin{align}
&\frac{1}{2\pi}\int_0^l\frac{(g'(s)+\overline{g'(s)})ds}{s-s_0}-\frac{(\kappa-1)}{2\pi i}\int_0^l\frac{(q(s)-\overline{q(s)})ds}{s-s_0}
\notag\\
&+\frac{\gamma_1}{4\mu}\kappa_0(s_0)\left\{\kappa_0(s_0)\left[\frac{\kappa-1}{2\pi}\int_0^l\frac{(g'(s)+\overline{g'(s)})ds}{s-s_0}+\frac{2\kappa}{\pi i}\int_0^l\frac{(q(s)-\overline{q(s)})ds}{s-s_0} \right]\right.
\notag\\
&+\left.i\frac{d}{ds_0}\left[\frac{\kappa-1}{2\pi}\int_0^l\frac{(g'(s)-\overline{g'(s)})ds}{s-s_0}+\frac{2\kappa}{\pi i}\int_0^l\frac{(q(s)+\overline{q(s)})ds}{s-s_0} \right]    \right\}
\label{3_10}\\
&=\R M_1(g',q)(s_0),\,\,\,\, s_0\in[0,l];
\notag\\
&\frac{1}{2\pi i}\int_0^l\frac{(g'(s)-\overline{g'(s)})ds}{s-s_0}+\frac{(\kappa-1)}{2\pi}\int_0^l\frac{(q(s)+\overline{q(s)})ds}{s-s_0}
\notag\\
&+\frac{\gamma_1}{4\mu}\frac{d}{ds_0}\left\{\kappa_0(s_0)\left[\frac{\kappa-1}{2\pi}\int_0^l\frac{(g'(s)+\overline{g'(s)})ds}{s-s_0}+\frac{2\kappa}{\pi i}\int_0^l\frac{(q(s)-\overline{q(s)})ds}{s-s_0} \right]\right.
\notag\\
&+\left.i\frac{d}{ds_0}\left[\frac{\kappa-1}{2\pi}\int_0^l\frac{(g'(s)-\overline{g'(s)})ds}{s-s_0}+\frac{2\kappa}{\pi i}\int_0^l\frac{(q(s)+\overline{q(s)})ds}{s-s_0} \right]    \right\}
\label{3_11}\\
&=\I M_1(g',q)(s_0),\,\,\,\, s_0\in[0,l];
\notag
\end{align}
where $M_1(g',q)(s_0)$ is a Fredholm integral operator of the first kind
\begin{align}
&M_1(g',q)(s_0)=-\frac{1}{2\pi}\int_0^l \left(k_1(s,s_0)g'(s)+k_2(s,s_0)\overline{g'(s)}-2ik_3(s,s_0)q(s)+2ik_2(s,s_0)\overline{q(s)}\right)ds
\notag\\
&-\frac{\gamma_1}{2\mu}(\kappa_0^2(s_0)+i\kappa'_0(s_0))\R\left [\frac{1}{2\pi}\int_0^l (k_4(s,s_0)g'(s)-k_2(s,s_0)\overline{g'(s)})ds\right.
\notag\\
&+\left.\frac{1}{\pi i}\int_0^l (\kappa k_1(s,s_0)q(s)+k_2(s,s_0)\overline{q(s)})ds \right]
\label{3_11a}\\
&-\frac{\gamma_1}{4\pi \mu}\kappa_0(s_0)\int_0^l \left(ik_{4,2}(s,s_0)g'(s)-ik_{2,2}(s,s_0)\overline{g'(s)}+2\kappa k_{1,2}(s,s_0)q(s)+2k_{2,2}(s,s_0)\overline{q(s)}\right)ds
\notag\\
&+i\frac{\gamma_1}{2\mu}\I\left[\frac{1}{2\pi}\int_0^l (k_{4,22}(s,s_0)g'(s)-k_{2,22}(s,s_0)\overline{g'(s)})ds \right.
\notag\\
&+\left.\frac{1}{\pi i}\int_0^l (\kappa k_{1,22}(s,s_0)q(s)+k_{2,22}(s,s_0)\overline{q(s)})ds \right]+(\kappa+1)f(s_0),
\notag
\end{align}
where $k_{j,2}(s,s_0)$ and $k_{j,22}(s,s_0)$ denote the first and the second partial derivatives of the kernels (\ref{3_8}) by the variable $s_0$. Here and further, assume for simplicity that $\kappa_0(s_0)\neq 0$ everywhere on $L$.

Similarly, by substituting (\ref{3_4}) into (\ref{3_9}) we obtain the following relations between the functions $g'(s)$ and $q(s)$:
\begin{align}
q(s_0)+\overline{q(s_0)}&=-\frac{i\gamma_1}{4\mu}\kappa_0(s_0)\left[\kappa_0(s_0)(g'(s_0)-\overline{g'(s_0)})+i(g''(s_0)+\overline{g''(s_0)}) \right];
\label{3_12}\\
q(s_0)-\overline{q(s_0)}&=\frac{\gamma_1}{4\mu}\frac{d}{ds_0}\left[\kappa_0(s_0)(g'(s_0)-\overline{g'(s_0)})+i(g''(s_0)+\overline{g''(s_0)})\right]
\notag\\
&=i\frac{d}{ds_0}\left[\frac{1}{\kappa_0(s_0)}(q(s_0)+\overline{q_0(s_0)}) \right].
\notag
\end{align}

Observe that the expression in the figure brackets $\{\ldots\}$ is the same in both of the equations (\ref{3_10}), (\ref{3_11}). Solve the equation (\ref{3_10}) for this expression and substitute into the equation (\ref{3_11}). Then, after differentiating and inverting the singular integrals, we arrive at the following regularized system of integro-differential equations:
\begin{align}
&\left(\frac{1}{2\kappa_0(s_0)}+\frac{\gamma_1(\kappa-1)}{8\mu}\kappa_0(s_0)\right)(g'(s_0)+\overline{g'(s_0)})+i\left(\frac{\kappa-1}{2\kappa_0(s_0)}-\frac{\gamma_1\kappa}{2\mu}\kappa_0(s_0)\right)(q(s_0)-\overline{q(s_0)})
\notag\\
&+\frac{\gamma_1}{4\mu}\left[i\frac{\kappa-1}{2}(g''(s_0)-\overline{g''(s_0)})+2\kappa(q'(s_0)+\overline{q'(s_0)}) \right]
\notag\\
&=-\frac{1}{\sqrt{s_0(l-s_0)}}\frac1{\pi}\int_0^l\left(\frac1{\kappa_0(s)}\R M_1(g',q)(s)
+M_2(g',q)(s) \right)\frac{\sqrt{s(l-s)}ds}{s-s_0}
\notag\\
&+\frac{D_1}{\pi \sqrt{s_0(l-s_0)}},\,\,\,\,s_0\in[0,l];
\label{3_13}\\
&-i\frac{\kappa_0(s_0)}{2}(g'(s_0)-\overline{g'(s_0)})+\frac{\kappa-1}{2}\kappa_0(s_0)(q(s_0)+\overline{q(s_0)})
\notag\\
&+\frac{\kappa'_0(s_0)}{2\kappa_0(s_0)}(g'(s_0)+\overline{g'(s_0)})+\frac{i(\kappa-1)}{2}\frac{\kappa'_0(s_0)}{\kappa_0(s_0)}(q(s_0)-\overline{q(s_0)})
\notag\\
&-\frac12(g''(s_0)+\overline{g''(s_0)})-\frac{i(\kappa-1)}{2}(q'(s_0)-\overline{q'(s_0)})
\notag\\
&=-\left.\frac{1}{\sqrt{s_0(l-s_0)}}\frac1{\pi}\int_0^l\right[\kappa_0(s)\I M_1(g',q)(s)-\kappa_0(s)\frac{d}{ds}\left(\frac{1}{\kappa_0(s)}\R M_1(g',q)(s)\right)
\notag\\
&+M_3(g',q) \left]\frac{\sqrt{s(l-s)}ds}{s-s_0}\right.+\frac{D_2}{\pi\sqrt{s_0(l-s_0)}},\,\,\,\,s_0\in[0,l];
\label{3_14}
\end{align}
where the regular kernels $M_2(g',q)(s)$ and $M_3(g',q)(s)$ have the following form:
\begin{align*}
M_2(g',q)(s_0)&=\frac{1}{2\pi}\int_0^l\left(\frac{1}{\kappa_0(s)}-\frac{1}{\kappa_0(s_0)} \right)\frac{(g'(s)+\overline{g'(s)})ds}{s-s_0}\\
&+\frac{i(\kappa-1)}{2\pi}\int_0^l\left(\frac{1}{\kappa_0(s)}-\frac{1}{\kappa_0(s_0)} \right)\frac{(q(s)-\overline{q(s)})ds}{s-s_0}\\
&+\frac{\gamma_1}{4\mu}\left[\frac{\kappa-1}{2\pi}\int_0^l(\kappa^2_0(s)-\kappa^2_0(s_0))\frac{(g'(s)+\overline{g'(s)})ds}{s-s_0}\right.\\
&-\left.\frac{2i\kappa}{\pi}\int_0^l(\kappa_0(s)-\kappa_0(s_0))\frac{(q(s)-\overline{q(s)})ds}{s-s_0}\right];\\
\end{align*}
\begin{align*}
M_3(g',q)(s_0)&=\frac{1}{2\pi i}\int_0^l(\kappa_0(s)-\kappa_0(s_0))\frac{(g'(s)-\overline{g'(s)})ds}{s-s_0}\\
&+\frac{(\kappa-1)}{2\pi}\int_0^l(\kappa_0(s)-\kappa_0(s_0))\frac{(q(s)+\overline{q(s)})ds}{s-s_0}\\
&+\frac{1}{2\pi}\int_0^l\left(\frac{\kappa'_0(s)}{\kappa_0(s)}-\frac{\kappa'_0(s_0)}{\kappa_0(s_0)}\right)\frac{(g'(s)+\overline{g'(s)})ds}{s-s_0}\\
&+\frac{i(\kappa-1)}{2\pi}\int_0^l\left(\frac{\kappa'_0(s)}{\kappa_0(s)}-\frac{\kappa'_0(s_0)}{\kappa_0(s_0)}\right)\frac{(q(s)-\overline{q(s)})ds}{s-s_0}.
\end{align*}
The constants $D_1$ and $D_2$ in the equations (\ref{3_13}) and (\ref{3_14}) can be specified by integrating these equations and using the additional condition (\ref{3_9a}): 
\begin{align}
D_1&=\int_0^l\left\{\frac{1}{2\kappa_0(s_0)}(g'(s_0)+\overline{g'(s_0)})+\frac{i(\kappa-1)}{2\kappa_0(s_0)}(q(s_0)-\overline{q(s_0)})\right.
\notag\\
&+\frac{\gamma_1}{4\mu}\kappa_0(s_0)\left[\frac{\kappa-1}{2}(g'(s_0)+\overline{g'(s_0)})-2i\kappa(q(s_0)-\overline{q(s_0)}) \right]
\label{3_14a}\\
&+\left.\frac{\gamma_1}{4\mu}\left[i\frac{\kappa-1}{2}(g''(s_0)-\overline{g''(s_0)})+2\kappa(q'(s_0)+\overline{q'(s_0)}) \right]\right\}ds_0
+\R\int_0^lg'(s_0)t'(s_0)ds_0;
\notag\\
D_2&=\int_0^l\left\{-i\frac{\kappa_0(s_0)}{2}(g'(s_0)-\overline{g'(s_0)})+\frac{\kappa-1}{2}\kappa_0(s_0)(q(s_0)+\overline{q(s_0)})\right.
\notag\\
&+\frac{\kappa'_0(s_0)}{2\kappa_0(s_0)}(g'(s_0)+\overline{g'(s_0)})+\frac{i(\kappa-1)}{2}\frac{\kappa'_0(s_0)}{\kappa_0(s_0)}(q(s_0)-\overline{q(s_0)})
\notag\\
&-\left.(g''(s_0)+\overline{g''(s_0)})-\frac{i(\kappa-1)}{2}(q'(s_0)-\overline{q'(s_0)})\right\}ds_0+\I\int_0^lg'(s_0)t'(s_0)ds_0.
\end{align}
Taking the constants $D_1$ and $D_2$ in the form (\ref{3_14a}) allows to satisfy the conditions (\ref{3_9a}) automatically.
Note that the system (\ref{3_13}) and (\ref{3_14}) has a H\"older continuous solution $g'(s)$, $q(s)$ only if the following conditions are satisfied at the ends of the crack:
\begin{align}
\left.\left(i\frac{\kappa-1}{2}(g'(s)-\overline{g'(s)})+2\kappa(q(s)+\overline{q(s)})\right)\right|_{s=0,l}&=0;
\label{3_15}\\
\left.\left(\frac{1}{2}(g'(s)+\overline{g'(s)})+\frac{i(\kappa-1)}{2}(q(s)-\overline{q(s)})\right)\right|_{s=0,l}&=0.
\notag
\end{align}
The physical meaning of the conditions (\ref{3_15}) is that the tensile stress $\sigma_n$ and the derivative of the displacement $du_2/ds$ are bounded at the tips of the crack.

Finally, we will reduce the system of two regularized integro-differential equations (\ref{3_13}), (\ref{3_14}) to the system of two weakly singular Fredholm integral equations. In the spirit of \cite{MikhPros1986}, denote the highest order derivatives in the equations (\ref{3_13}), (\ref{3_14}) as new unknown functions:
\begin{equation}
\label{3_16}
\varphi(s_0)=i(g''(s_0)-\overline{g''(s_0)});\,\,\,\,\,\,\psi(s_0)=q''(s_0)+\overline{q''(s_0)}.
\end{equation}
Integrating the formulas (\ref{3_16}) and making use of the formulas (\ref{3_12}) obtain:
\begin{align}
&q'(s_0)+\overline{q'(s_0)}=\int_0^l\omega_0(s,s_0)\psi(s)ds+C_1;
\notag\\
&q(s_0)+\overline{q(s_0)}=\int_0^l\omega_1(s,s_0)\psi(s)ds+C_1s_0+C_2;
\notag\\
&i(g'(s_0)-\overline{g'(s_0)})=\int_0^l\omega_0(s,s_0)\varphi(s)ds+C_3;
\notag\\
&g''(s_0)+\overline{g''(s_0)}=\frac{4\mu}{\gamma_1}\frac{1}{\kappa_0(s_0)}\left(\int_0^l\omega_1(s,s_0)\psi(s)ds+C_1s_0+C_2 \right)
\notag\\
&+\kappa_0(s_0)\left(\int_0^l\omega_0(s,s_0)\varphi(s)ds+C_3 \right);
\label{3_17}\\
&g'(s_0)+\overline{g'(s_0)}=\frac{4\mu}{\gamma_1}\int_0^l\frac{\omega_0(s,s_0)}{\kappa_0(s)}\left(\int_0^l\omega_1(s_1,s)\psi(s_1)ds_1+C_1s+C_2 \right)ds
\notag\\
&+\int_0^l\omega_0(s,s_0)\kappa_0(s)\left(\int_0^l\omega_0(s_1,s)\varphi(s_1)ds_1+C_3 \right)ds+C_4;
\notag\\
&i(q'(s_0)-\overline{q'(s_0)})=-\frac1{\kappa_0(s_0)}\psi(s_0)+\frac{2\kappa'_0(s_0)}{\kappa_0^2(s_0)}\left(\int_0^l\omega_0(s,s_0)\psi(s)ds+C_1 \right)
\notag\\
&+\frac{d}{ds_0}\left(\frac{\kappa'_0(s_0)}{\kappa_0^2(s_0)}\right)\left(\int_0^l\omega_1(s,s_0)\psi(s)ds+C_1s_0+C_2\right);
\notag\\
&i(q(s_0)-\overline{q(s_0)})=-\frac{1}{\kappa_0(s_0)}\left(\int_0^l\omega_0(s,s_0)\psi(s)ds+C_1 \right)
\notag\\
&+\frac{\kappa'_0(s_0)}{\kappa_0^2(s_0)}\left(\int_0^l\omega_1(s,s_0)\psi(s)ds+C_1s_0+C_2\right),
\notag
\end{align}
where
$$
\omega_0(s,s_0)=\left\{
\begin{array}{cc}
1, & \mbox{ if } s\in[0,s_0],\\ 
0, & \mbox{ if } s\notin[0,s_0], 
\end{array}
\right.,\,\,\,
\omega_1(s,s_0)=\int_0^l\omega_0(s,s_1)\omega_0(s_1,s_0)ds_1.
$$

The four constants $C_1$, $C_2$, $C_3$ and $C_4$ can be found from the four additional conditions (\ref{3_15}). Substituting the formulas (\ref{3_17}) into the equations (\ref{3_13}), (\ref{3_14}) we obtain a system of two weakly singular Fredholm equations of the second kind with respect to the functions $\varphi(s)$, $\psi(s)$. This system has a unique solution with the exception of the discrete spectrum of the Fredholm operators. It also follows that the functions $\varphi(s)$ and $\psi(s)$ can have at most integrable power singularities of the order $1/2$ at the crack tips. Hence, the functions $g'(s)$ and $q(s)$ are bounded at those points. Thus, the stresses and the displacements can have at most logarithmic singularities at the crack tips. On the other hand, the conditions (\ref{3_15}) provide that the stress $\sigma_n^{\pm}$ and the derivative of the displacement $du_2^{\pm}/ds$ are bounded at the end points of the crack. Thus, only the functions $\tau_n^{\pm}$ and $du_1^{\pm}/ds$ may possess logarithmic singularities at the crack tips.

\section{Straight Crack}
While zero curvatures $\kappa_0(s_0)=0$ have been explicitly excluded from the consideration until now, it is possible to study a straight crack $\kappa_0(s_0)=0$ using the same methods as above. The equations (\ref{3_10}), (\ref{3_11}) in this particular case become:
\begin{align}
&\frac{1}{2\pi}\int_0^l\frac{(g'(s)+\overline{g'(s)})ds}{s-s_0}+\frac{i(\kappa-1)}{2\pi}\int_0^l\frac{(q(s)-\overline{q(s)})ds}{s-s_0}
\notag\\
&=\R M_1(g',q)(s_0),\,\,\,\, s_0\in[0,l];
\label{4_1}\\
&\frac{1}{2\pi i}\int_0^l\frac{(g'(s)-\overline{g'(s)})ds}{s-s_0}-\frac{\gamma_1(\kappa-1)}{8\pi i\mu}\int_0^l\frac{(g'''(s)-\overline{g'''(s)})ds}{s-s_0}=\I M_1(g',q)(s_0),\,\,\,\, s_0\in[0,l];
\notag
\end{align}
where the expression for the regular integral operator $M_1(g',q)(s_0)$ can be obtained from the formula (\ref{3_11a}) by substituting $\kappa_0(s_0)=0$. Furthermore, instead of the equations (\ref{3_12}) we obtain:
\begin{align}
q(s_0)+\overline{q(s_0)}&=0,\,\,\,\,s_0\in[0,l];
\label{4_2}\\
i(q(s_0)-\overline{q(s_0)})&=-\frac{\gamma_1}{4\mu}(g'''(s_0)+\overline{g'''(s_0)}),\,\,\,\,s_0\in[0,l].
\notag
\end{align}
The conditions (\ref{3_15}) become in this particular case:
\begin{align}
&g'(0)-\overline{g'(0)}=0,\,\,\,\,g''(0)-\overline{g''(0)}=0;
\label{4_3}\\
&g'(l)-\overline{g'(l)}=0,\,\,\,\,g''(l)-\overline{g''(l)}=0.
\notag
\end{align}

Just like in the previous section, it is possible to invert the singular integrals in the equations (\ref{4_1}) to obtain the system of two regular integro-differential equations. Analogously, introduce two new unknown functions:
\begin{equation}
\label{4_4}
\varphi(s_0)=g'''(s_0)+\overline{g'''(s_0)},\,\,\,\,\psi(s_0)=i(g'''(s_0)-\overline{g'''(s_0)}),\,\,\,\,s_0\in[0,l],
\end{equation}
and integrate them similarly to (\ref{3_17}). Substituting (\ref{4_4}) into the system of integro-differential equations (\ref{4_1}) we obtain a regularized system of two weakly singular Fredholm equations of the second kind. Similarly, it follows that $g'(s)$ and $q(s)$ are bounded at the crack tips. Furthermore, both the normal and the shear stresses $\sigma_n^{\pm}$ and $\tau_n^{\pm}$ and the derivatives of the displacements $du_2^{\pm}/ds$, $d^2u_2^{\pm}/ds^2$ are bounded at the crack tips, which completely corroborates the results of \cite{SendovaWalton2010} which were obtained there using different techniques.

\section{Numerical solution of the system (\ref{3_10}), (\ref{3_11})}

\begin{figure}[ht]
	\centering
		\scalebox{0.4}{\includegraphics{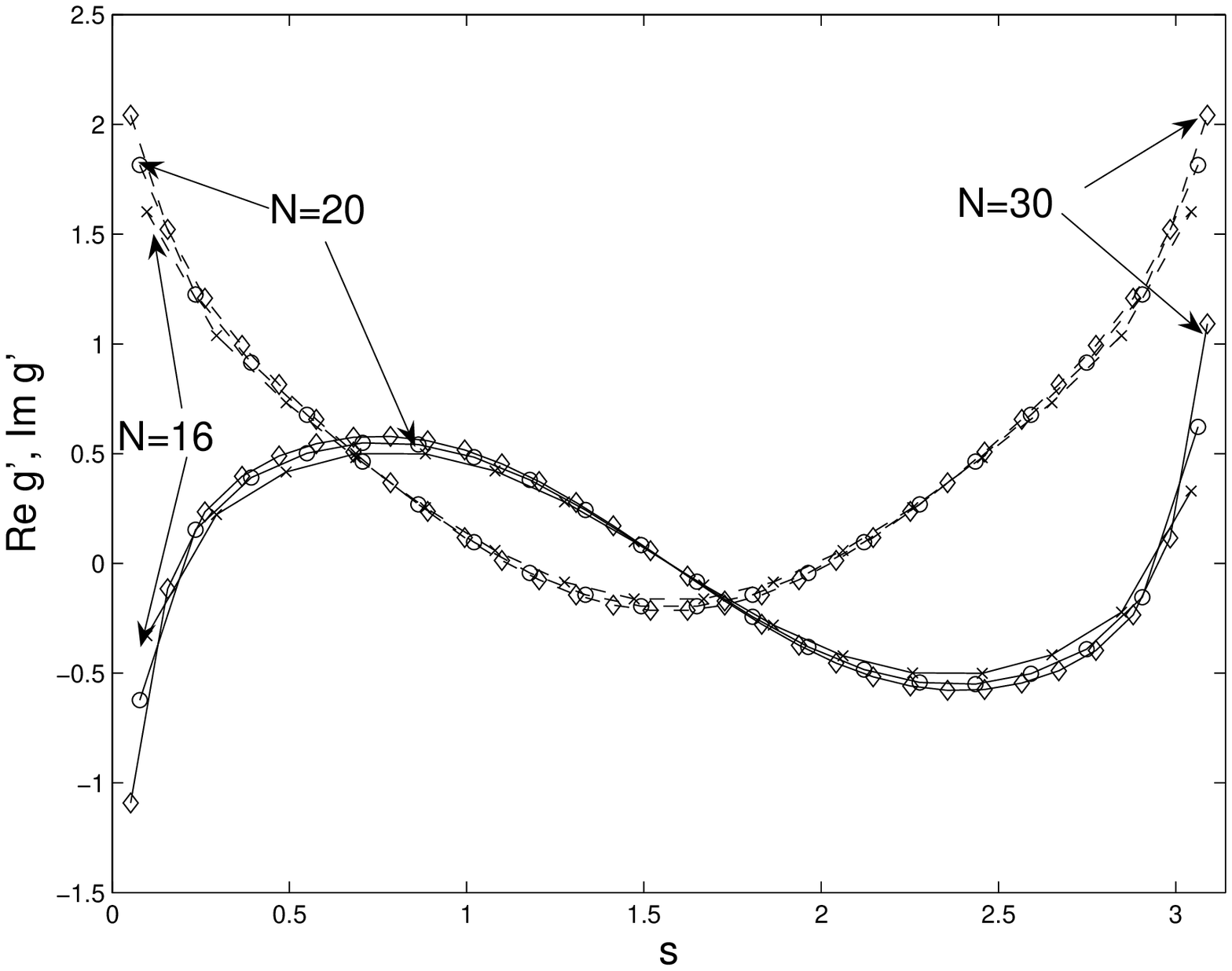}}
	\caption{Graphs of $\R g'$ and $\I g'$ for $N=16$, $N=20$ and $N=30$.}
	\label{fig1}
\end{figure}

\begin{figure}[ht]
	\centering
		\scalebox{0.35}{\includegraphics{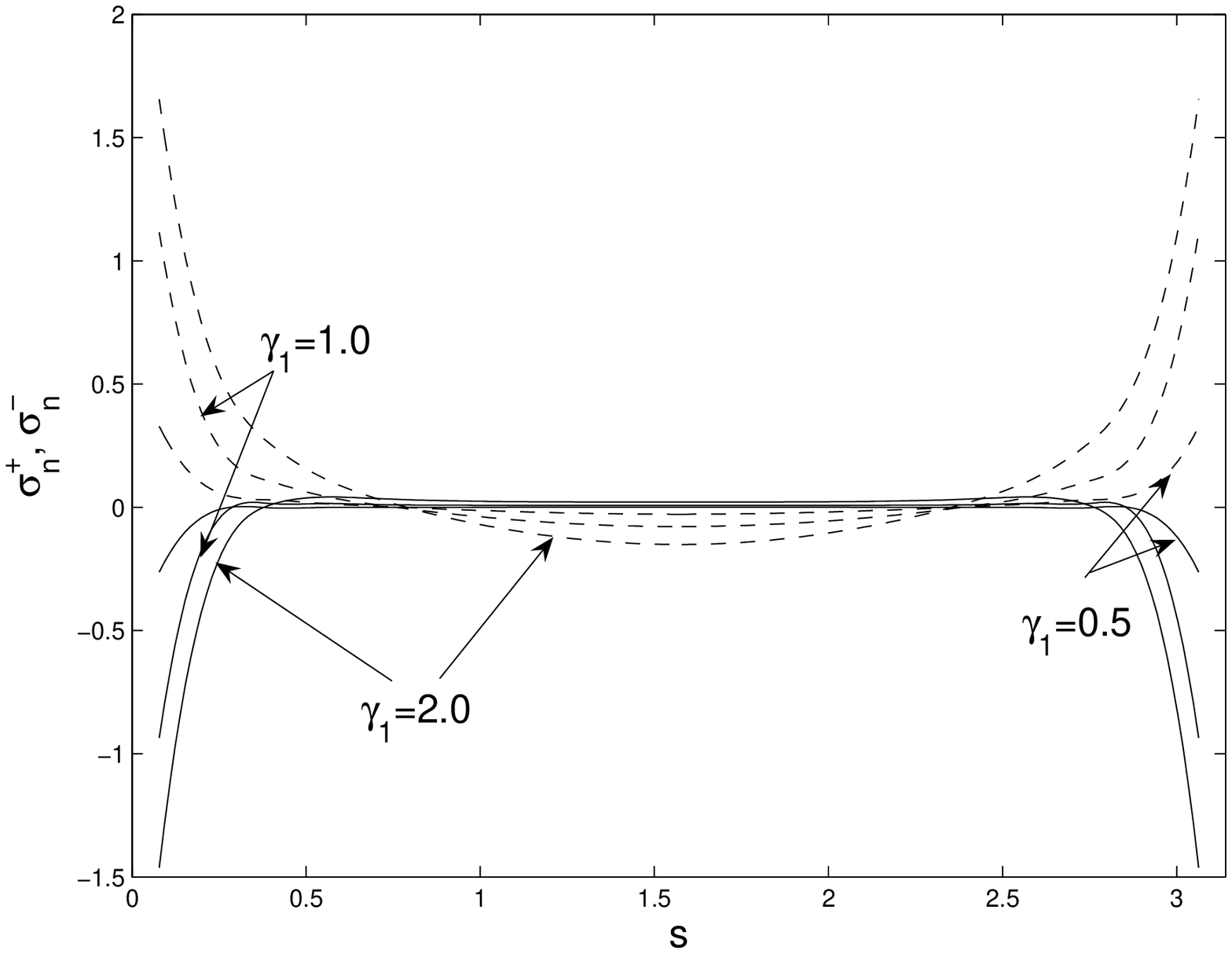} \includegraphics{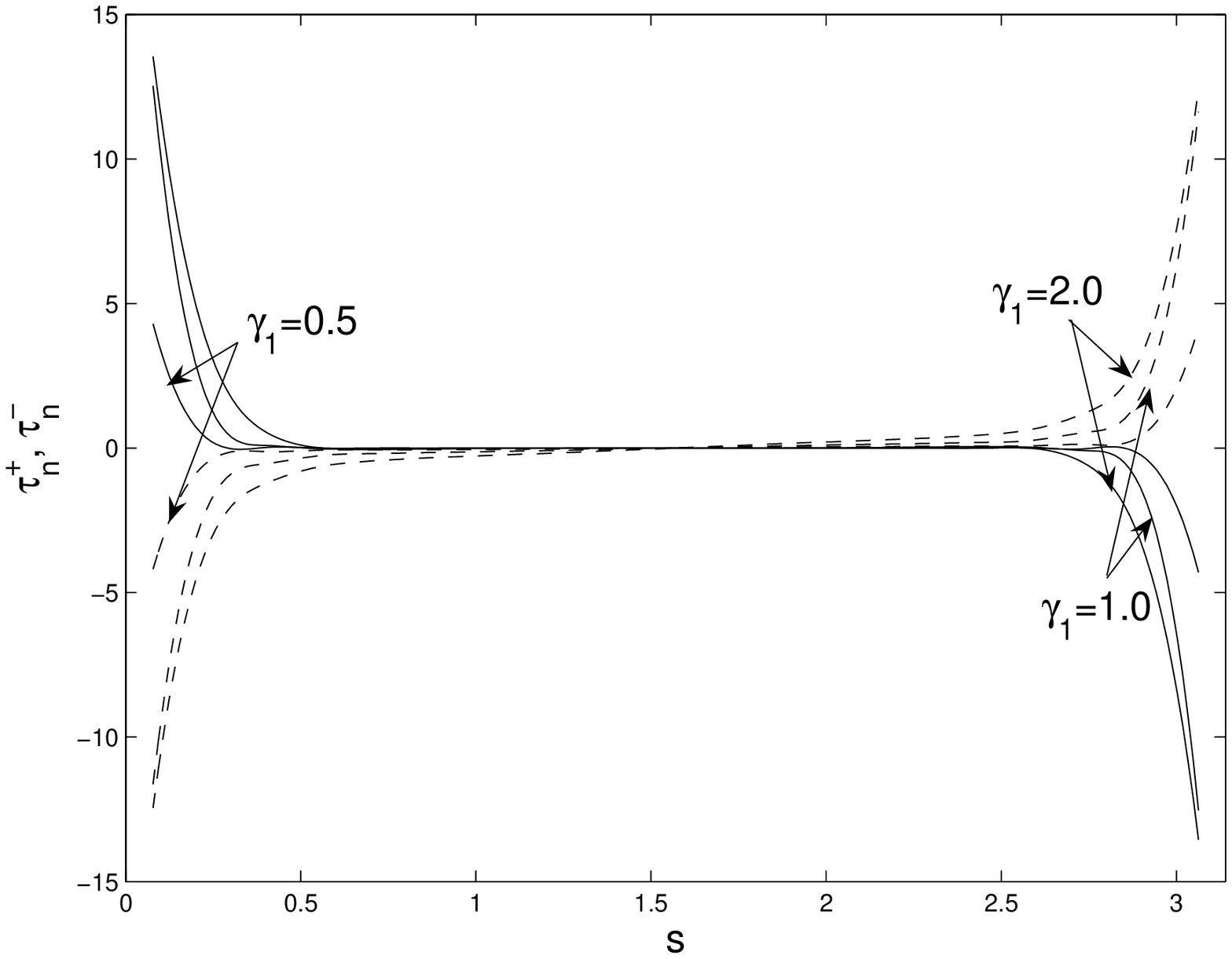}}
		\scalebox{0.35}{\includegraphics{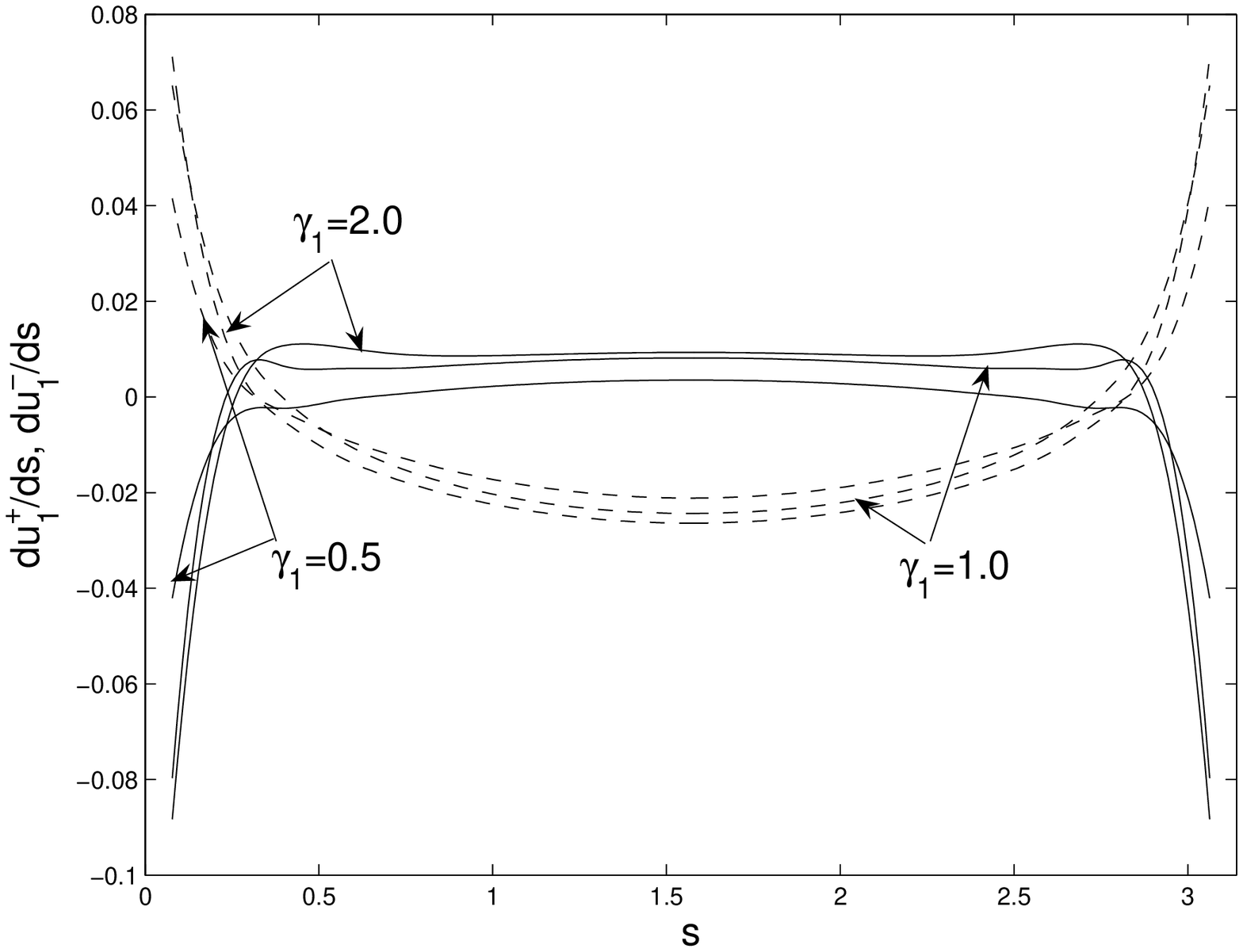} \includegraphics{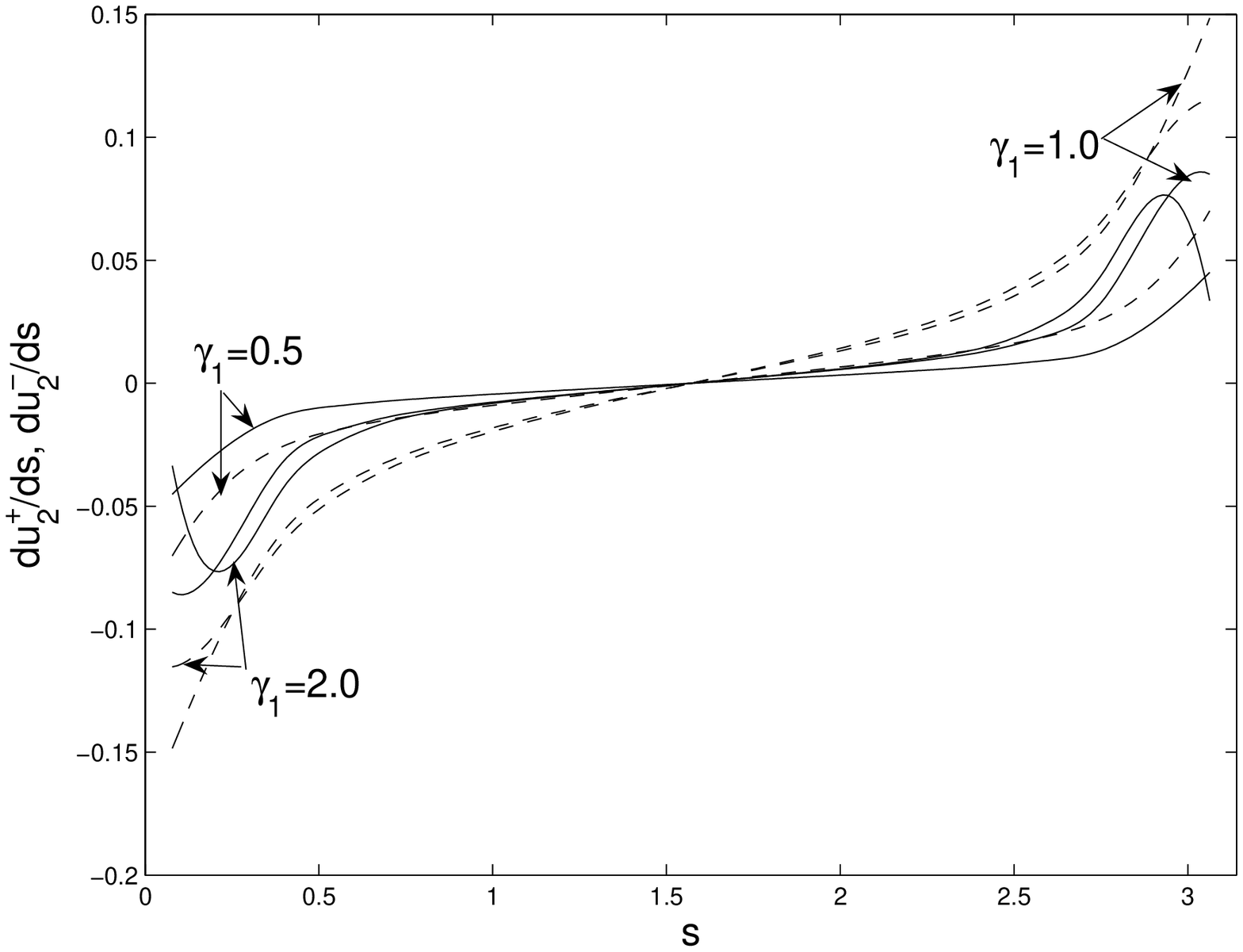}}
	\caption{Graphs of stresses $\sigma^{\pm}_n$, $\tau^{\pm}_n$ and derivatives of displacements $du_1^{\pm}/ds$, $du_2^{\pm}/ds$ for $\sigma_1^{\infty}=1$, $\sigma_2^{\infty}=0$.}
	\label{fig2}
\end{figure}

The system of two weakly singular Fredholm integral equations obtained by regularization of the system (\ref{3_13}), (\ref{3_14}) in general has a very complicated form and is inconvenient for numerical solution. Hence, for practical purposes it is easier to solve the initial system of integro-differential equations (\ref{3_10}), (\ref{3_11}). It is well known that the systems of singular integral and integro-differential equations of this kind are notoriously difficult to solve numerically with a good accuracy and a reasonable computational time. From the authors' point of view, one of the best ways to solve (\ref{3_10}), (\ref{3_11}) is to take a Taylor polynomial approximations of the unknown function $g'(s)$ in the following form:
\begin{equation}
\label{5_1}
g'(s)=\sum_{k=0}^{N}(g^1_k+ig^2_k)(s-l/2)^k,
\end{equation}
where $g^1_k$ and $g^2_k$ are real coefficients. Using the formulas (\ref{3_12}), one obtains the following expressions for the real and imaginary parts of the function $q(s)$:
\begin{align}
q(s)+\overline{q(s)}&=\frac{\gamma_1}{2\mu}\kappa_0(s)\sum_{k=0}^N\left[kg^1_k(s-l/2)^{k-1}+g_k^2\kappa_0(s)(s-l/2)^k \right],
\label{5_2}\\
i(q(s)-\overline{q(s)})&=-\frac{\gamma_1}{2\mu}\sum_{k=0}^N\left[k(k-1)g^1_k(s-l/2)^{k-2}+g_k^2(k\kappa_0(s)(s-l/2)^{k-1}\right.
&+\left.\kappa'_0(s)(s-l/2)^{k}) \right].
\notag
\end{align}
The following formula is used for numerical evaluation of the integrals in the equations (\ref{3_10}), (\ref{3_11}):
$$
\int_0^lK(s,s_0)\phi(s)ds=\frac{l}{N+1}\sum_{k=0}^NK(\tau_k,s_j)\phi(\tau_k),
$$
where $\tau_k=lk/N$, $k=0,1,\ldots,N$.

\begin{figure}[ht]
	\centering
		\scalebox{0.35}{\includegraphics{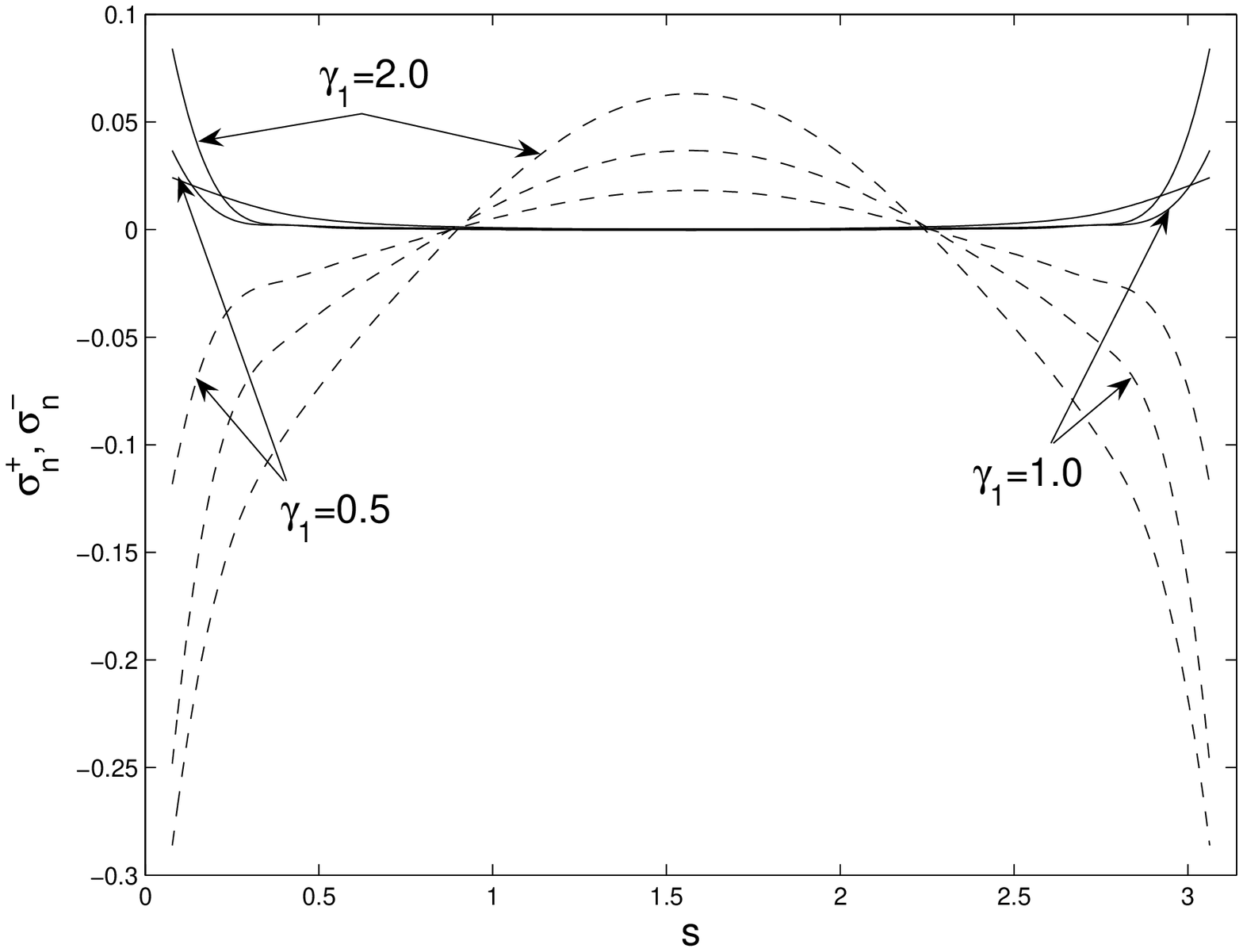} \includegraphics{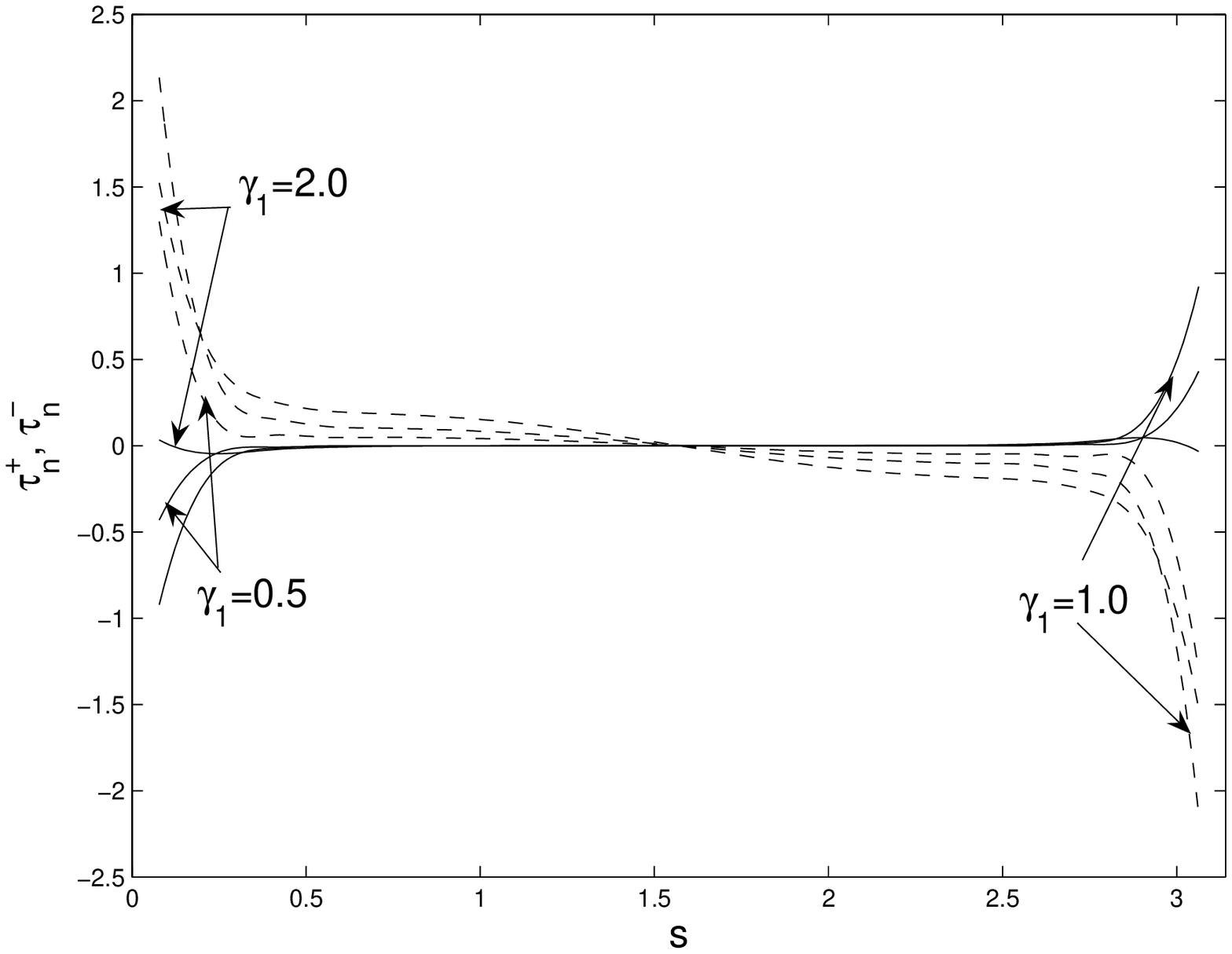}}
		\scalebox{0.35}{\includegraphics{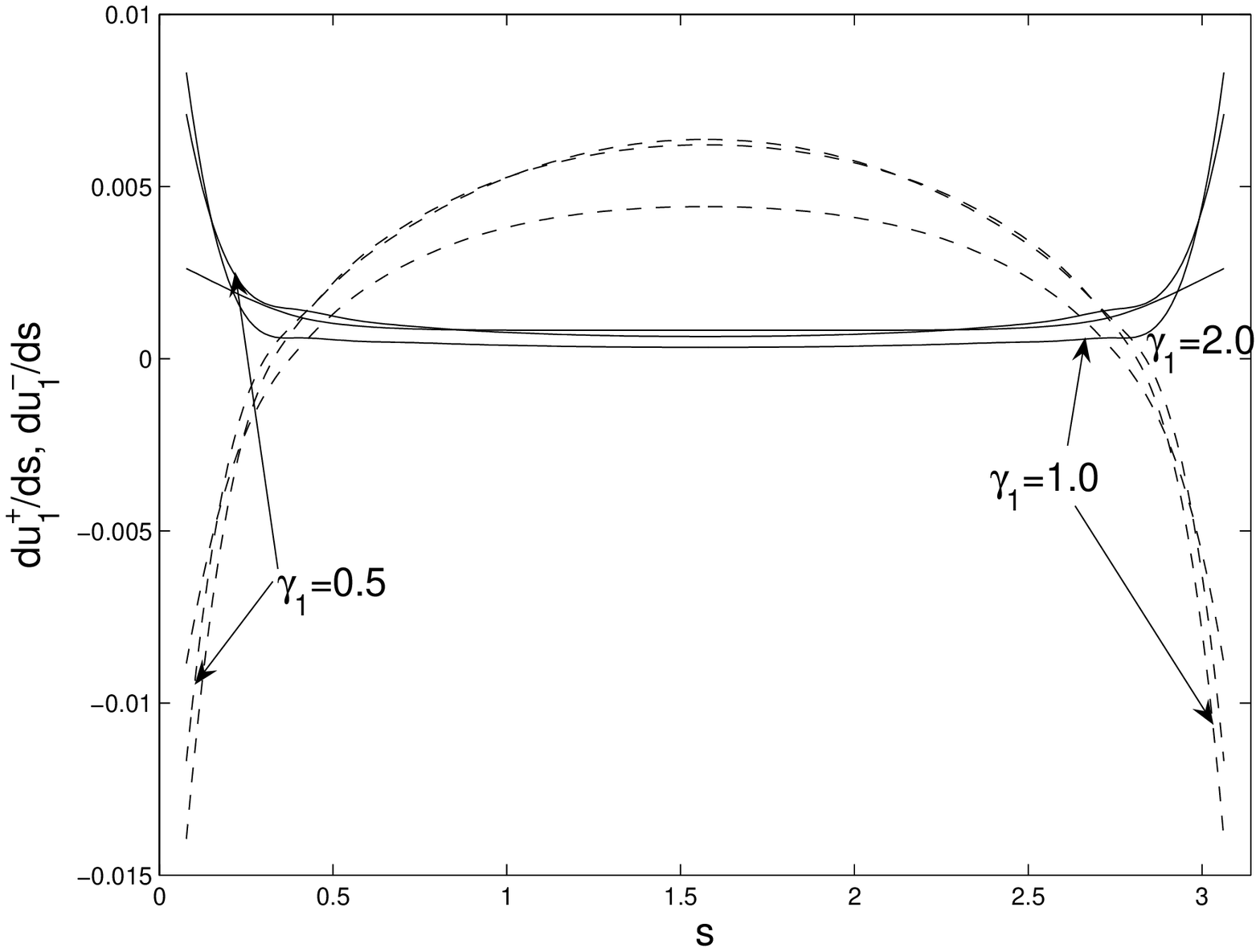} \includegraphics{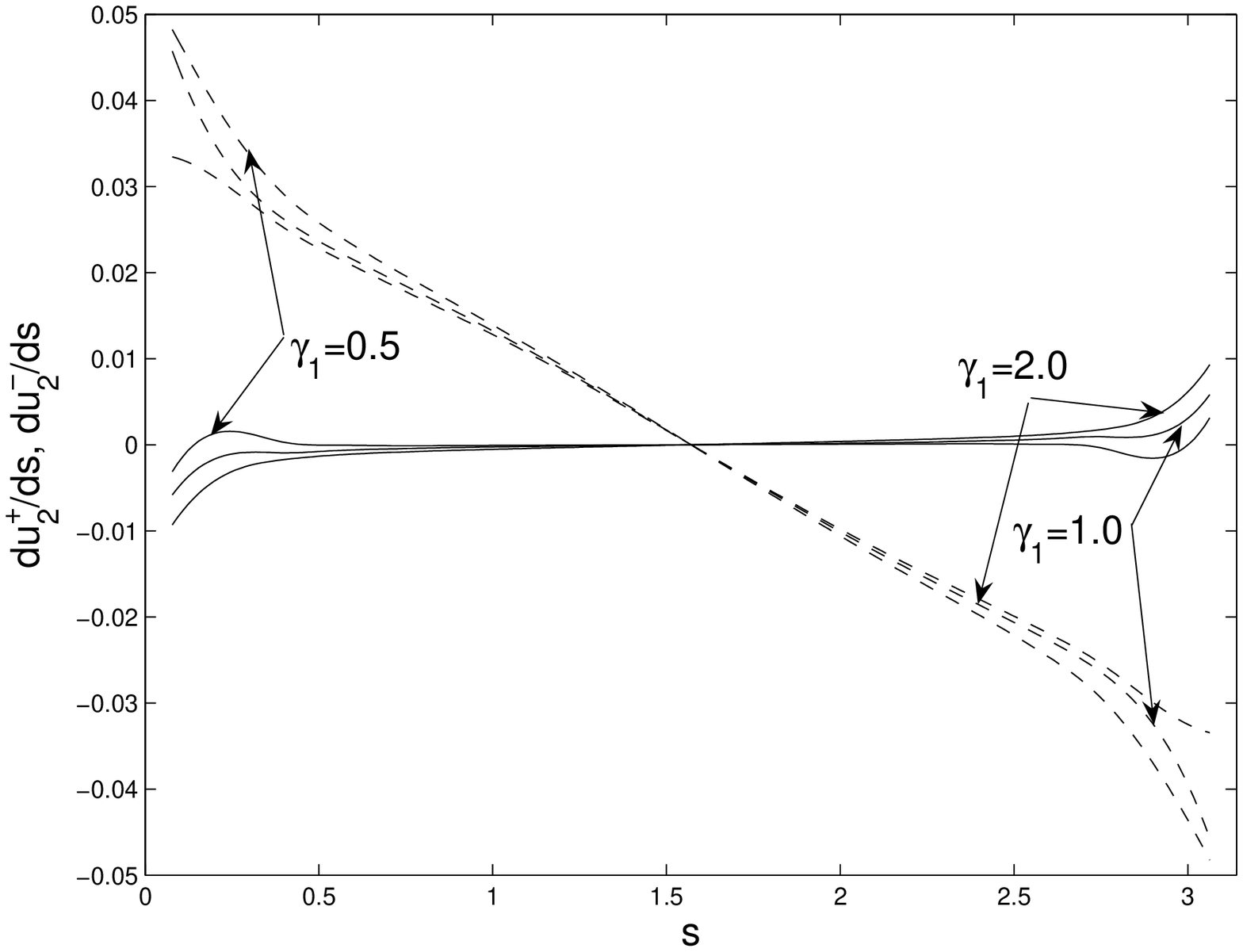}}
			\caption{Graphs of stresses $\sigma^{\pm}_n$, $\tau^{\pm}_n$ and derivatives of displacements $du_1^{\pm}/ds$, $du_2^{\pm}/ds$ for $\sigma_1^{\infty}=0$, $\sigma_2^{\infty}=1$.}
	\label{fig3}
\end{figure}

Substituting the representations (\ref{5_1}), (\ref{5_2}) into the singular integro-differential equations (\ref{3_10}), (\ref{3_11}) and the additional conditions (\ref{3_9a}) we obtain the following system of linear algebraic equations:
\begin{align}
&\sum_{k=0}^N (g^1_kN^1_k(s_j)+g^2_kN^2_k(s_j))=\R f(s_j),\,\,\,\,\,\sum_{k=0}^N (g^1_kN^3_k(s_j)+g^2_kN^4_k(s_j))=\I f(s_j),
\notag\\
&\sum_{k=0}^N(g^1_k+ig^2_k)\int_0^lt'(s)(s-l/2)^kds=0,
\label{5_3}
\end{align}
where $s_j=(2j-1)l/(2N)$, $j=1,2,\ldots,N$,  the coefficients $N^1_k(s_j)$, $N^2_k(s_j)$, $N^3_k(s_j)$ and $N^4_k(s_j)$ are obtained from the equations (\ref{3_10}), (\ref{3_11}). The solution of this system of linear algebraic equations provides the coefficients of the Taylor polynomials (\ref{5_1}).

\begin{figure}[ht]
	\centering
		\scalebox{0.35}{\includegraphics{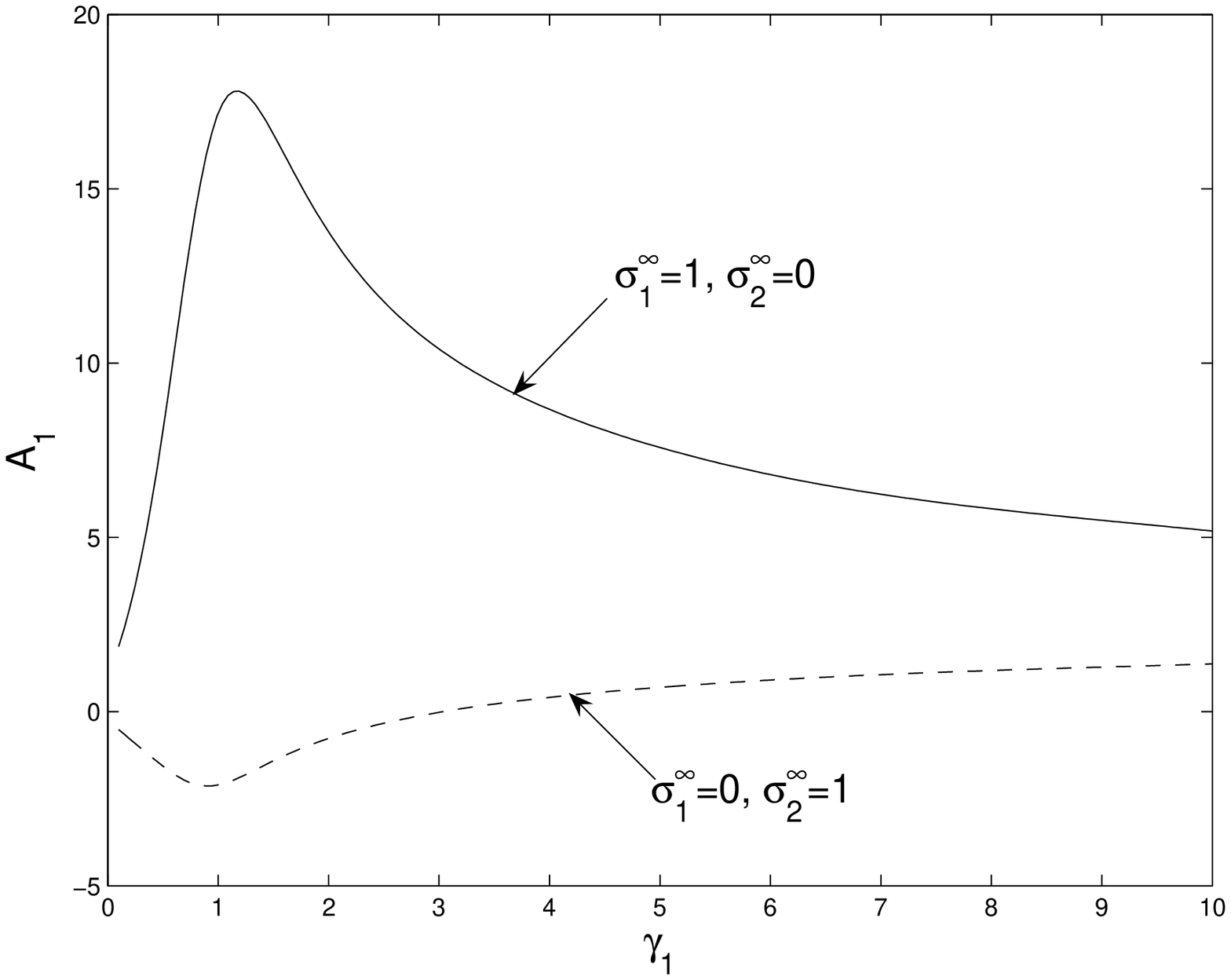}}
		\scalebox{0.35}{\includegraphics{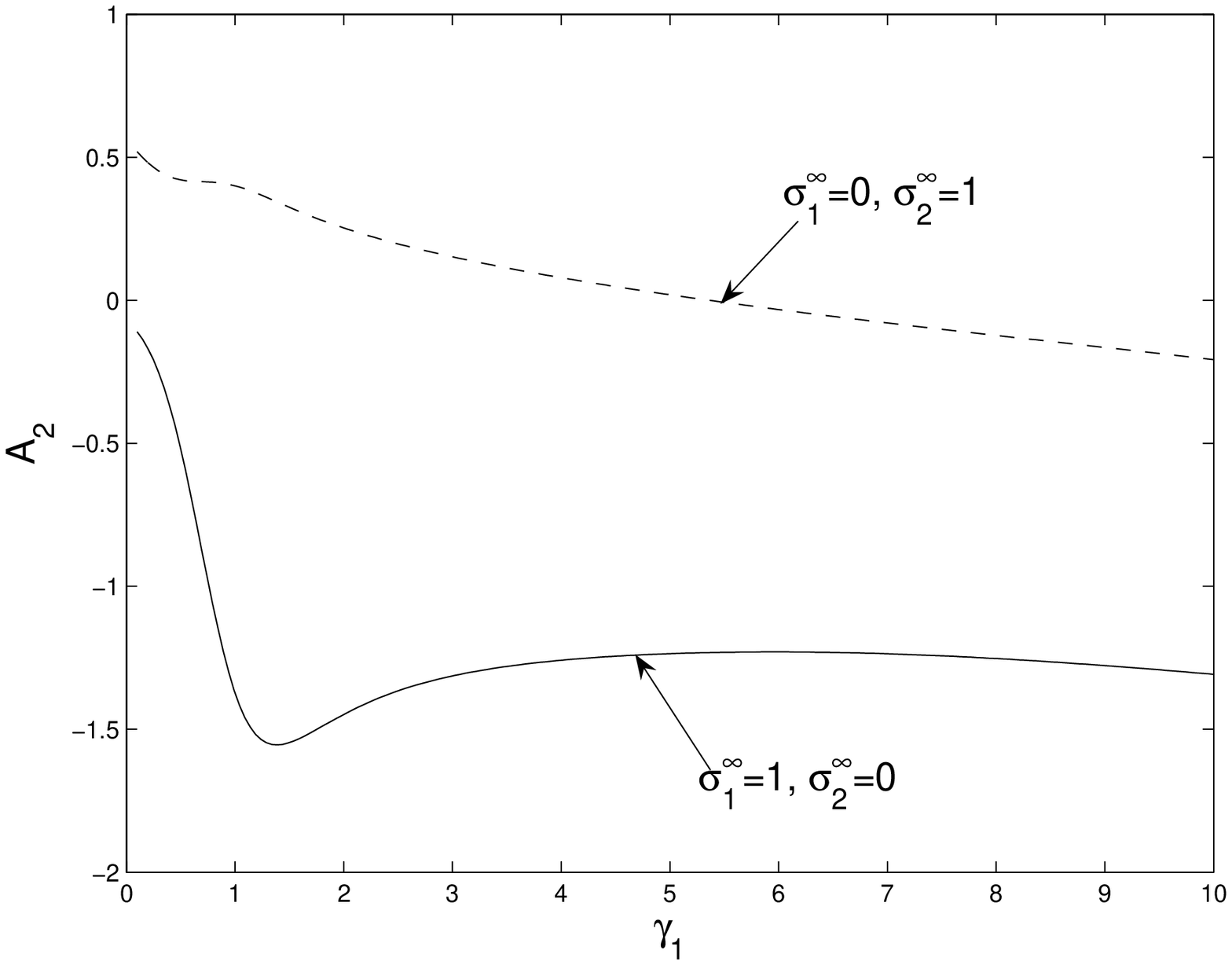}}
			\caption{Graphs of the dependence of the singularity coefficients $A_1$ and $A_2$ on the coefficient $\gamma_1$.}
	\label{fig4}
\end{figure}

Note that the system (\ref{3_10}), (\ref{3_11}) has additional conditions (\ref{3_15}) already built into it. It has been verified that the numerical solutions obtained from the system (\ref{5_3}) satisfies these conditions approximately. 

It has been also observed that this method provides good results even for relatively small values of $N$. For instance, fig. \ref{fig1} shows the graphs of the function $g'(s)$ in the case of the semicircular arc of a unit radius computed for $N=16$, $N=20$ and $N=30$ and the following values of the parameters $\mu=60$, $\kappa=2.5$, $\gamma_1=1.0$. It can be seen that a reasonable convergence is obtained already for $N=20$. Similar results have been observed in all considered cases. The proposed method of solution of the system (\ref{3_10}), (\ref{3_11}) is computationally efficient and leads to numerical solution of the relatively small systems of linear algebraic equations of the order $2N+2$, coefficients of which are computed by using numerical integration.

Observe also that since the function $g'(s)$ is bounded at the crack tips corresponding to the values of the parameter $s=0,l$, it follows from the equation (\ref{3_4}) that the crack has a sharp opening profile at the tips. This is consistent with experimental observation and is an improvement from the classical LEFM picture in which the crack opening profile is blunt. 

The proposed method of solution in general allows one to study an arbitrary smooth cracks. Some of the numerical results obtained using this method are presented below. 

\begin{figure}[ht]
	\centering
		\scalebox{0.35}{\includegraphics{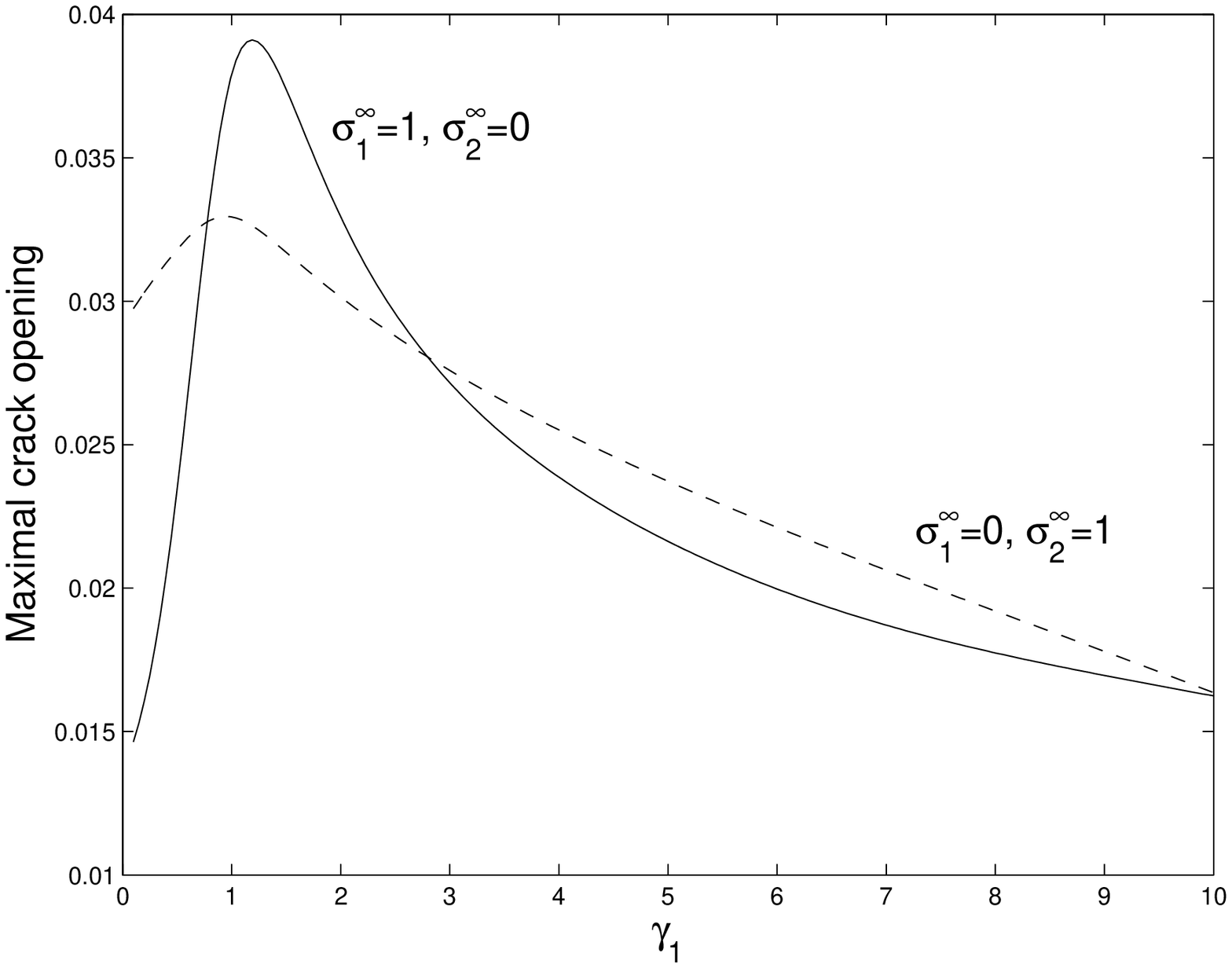}}
		\scalebox{0.35}{\includegraphics{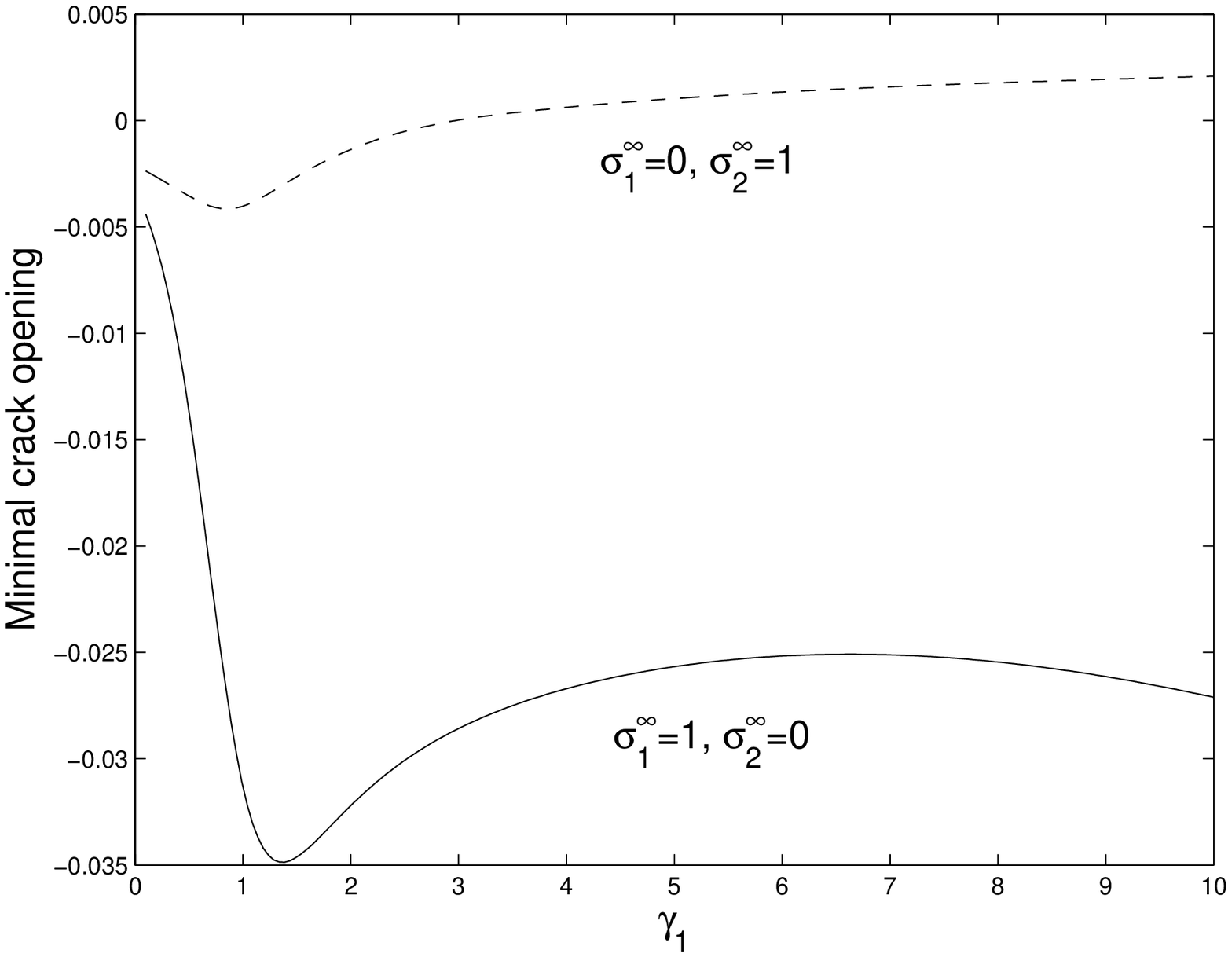}}
			\caption{Graphs of the dependence of the maximal and minimal crack opening on the coefficient $\gamma_1$.}
	\label{fig5}
\end{figure}

The first studied example is concerned with a semicircular crack of a unit radius parametrically defined by the equation $t(s)=e^{is}$, $s\in[0,\pi]$. We consider the following two types of the loading applied at the infinity: a horizontal ($\sigma_1^{\infty}=1$, $\sigma_2^{\infty}=0$) and a vertical stretching ($\sigma_1^{\infty}=0$, $\sigma_2^{\infty}=1$). The physical parameters of the material are $\mu=60$, $\kappa=2.5$. The dependence of the stresses and the derivatives of the displacements on the parameter $s$ for these two cases is shown on the figs. \ref{fig2} and \ref{fig3} correspondingly for three values of the proportionality constant $\gamma_1=0.5$, $\gamma_1=1.0$ and $\gamma_1=2.0$. Because of the symmetry of the crack geometry and the applied loading the stresses are distributed symmetrically with respect to the center of the crack. Observe that if the proportionality coefficient $\gamma_1=0$ then the conditions (\ref{3_9}) become $(\sigma_n+i\tau_n)^{\pm}(s_0)=0$ on the crack boundary. It can be seen from numerical results (and in particular figs. \ref{fig2} and \ref{fig3}) that $(\sigma_n+i\tau_n)^{\pm}(s_0) \to 0$ as $\gamma_1\to 0$.

\begin{figure}[ht]
	\centering
		\scalebox{0.4}{\includegraphics{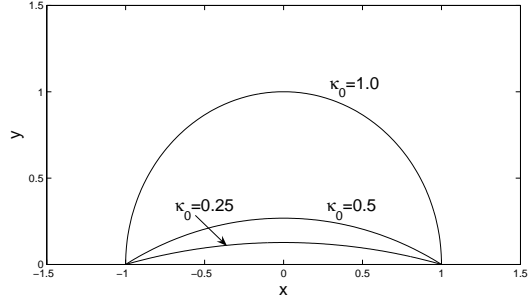}}
			\caption{Crack geometry for $\kappa_0=0.25$, $\kappa_0=0.5$ and $\kappa_0=1.0$.}
	\label{fig5a}
\end{figure}

Here the dashed lines correspond to the values of the stresses and displacements on the upper surface of the crack (denoted by ``$+$" sign), and the solid lines correspond to the values of the stresses and displacements on the lower surface of the crack (denoted by ``$-$" sign).

\begin{figure}[ht]
	\centering
		\scalebox{0.35}{\includegraphics{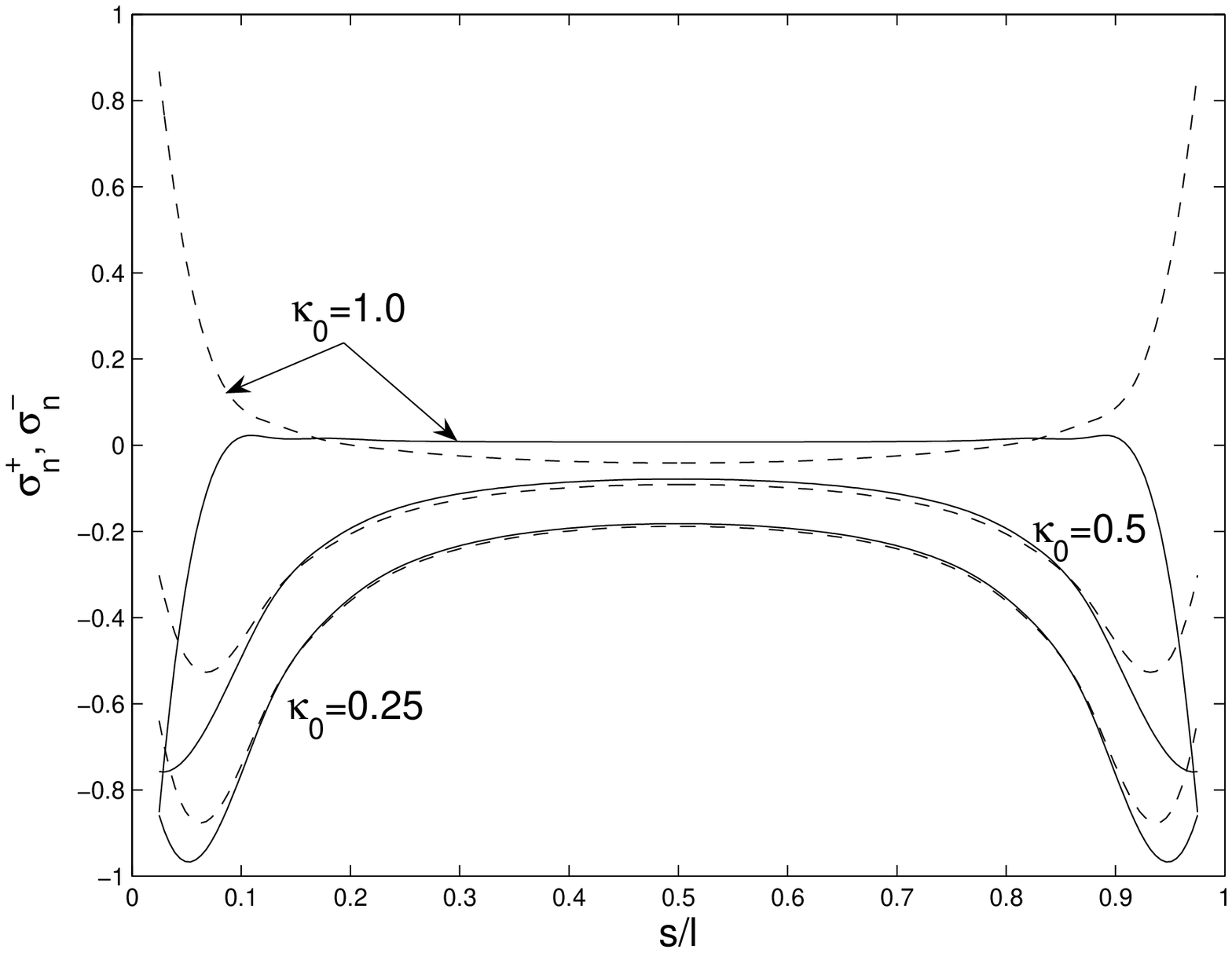} \includegraphics{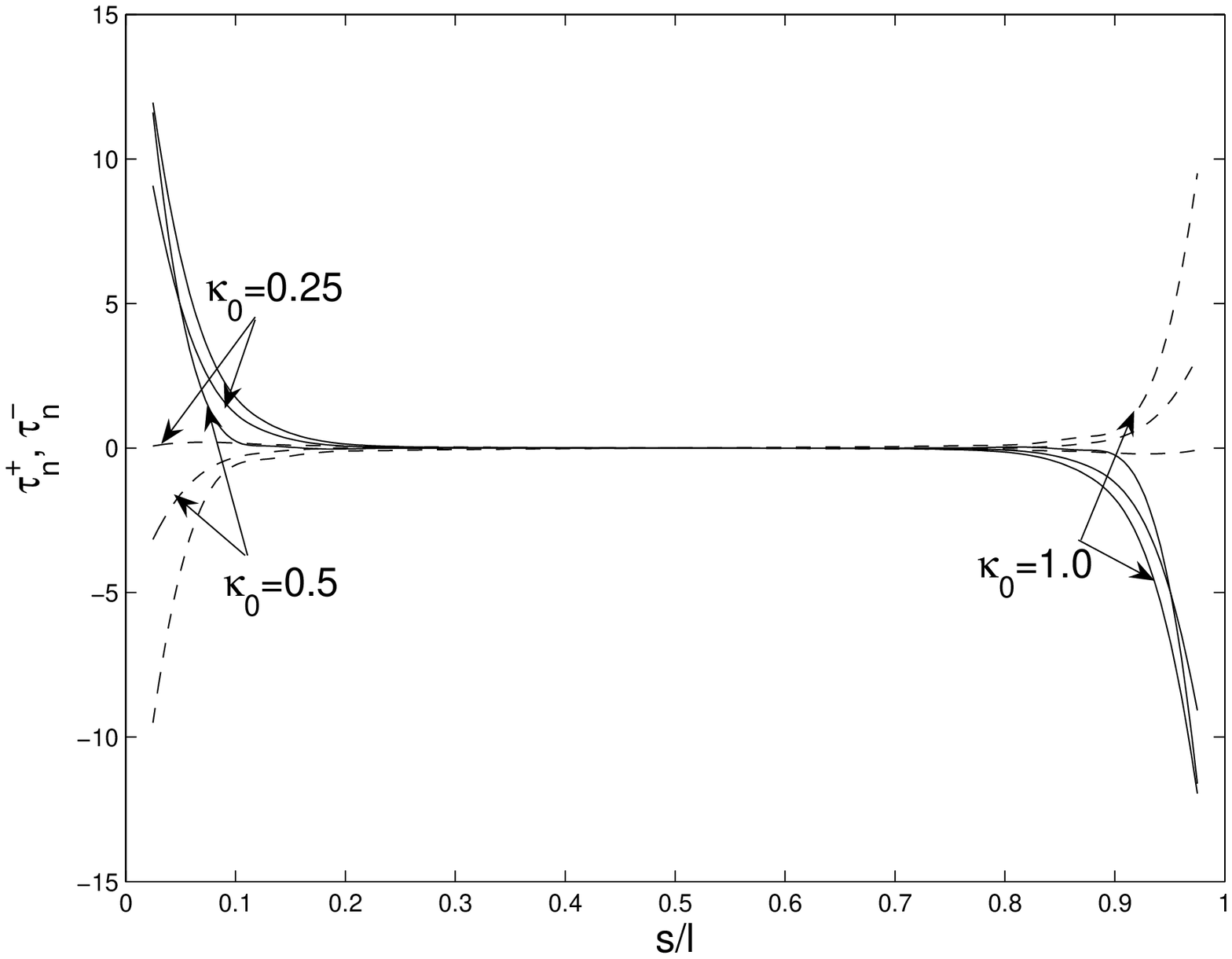}}
		\scalebox{0.35}{\includegraphics{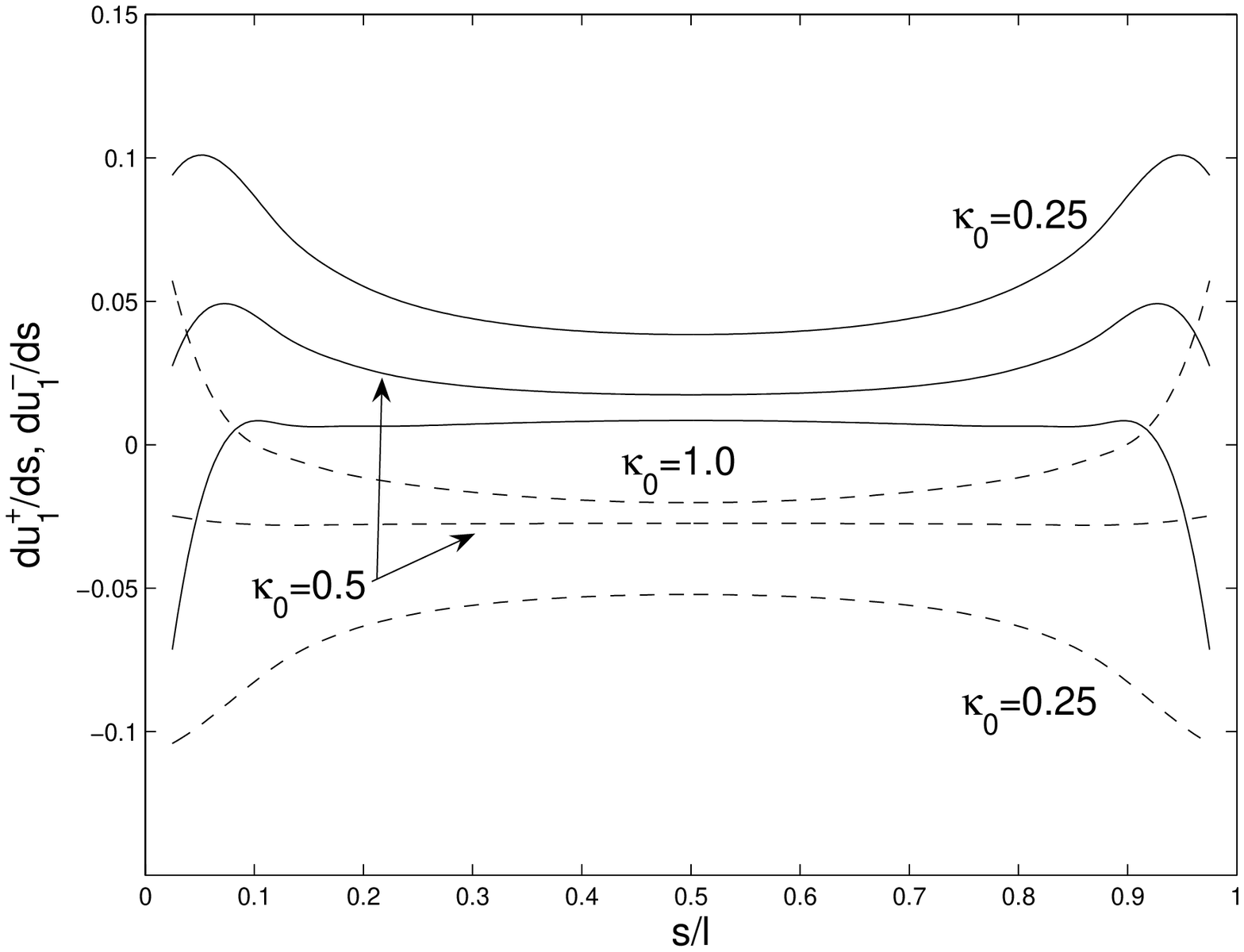} \includegraphics{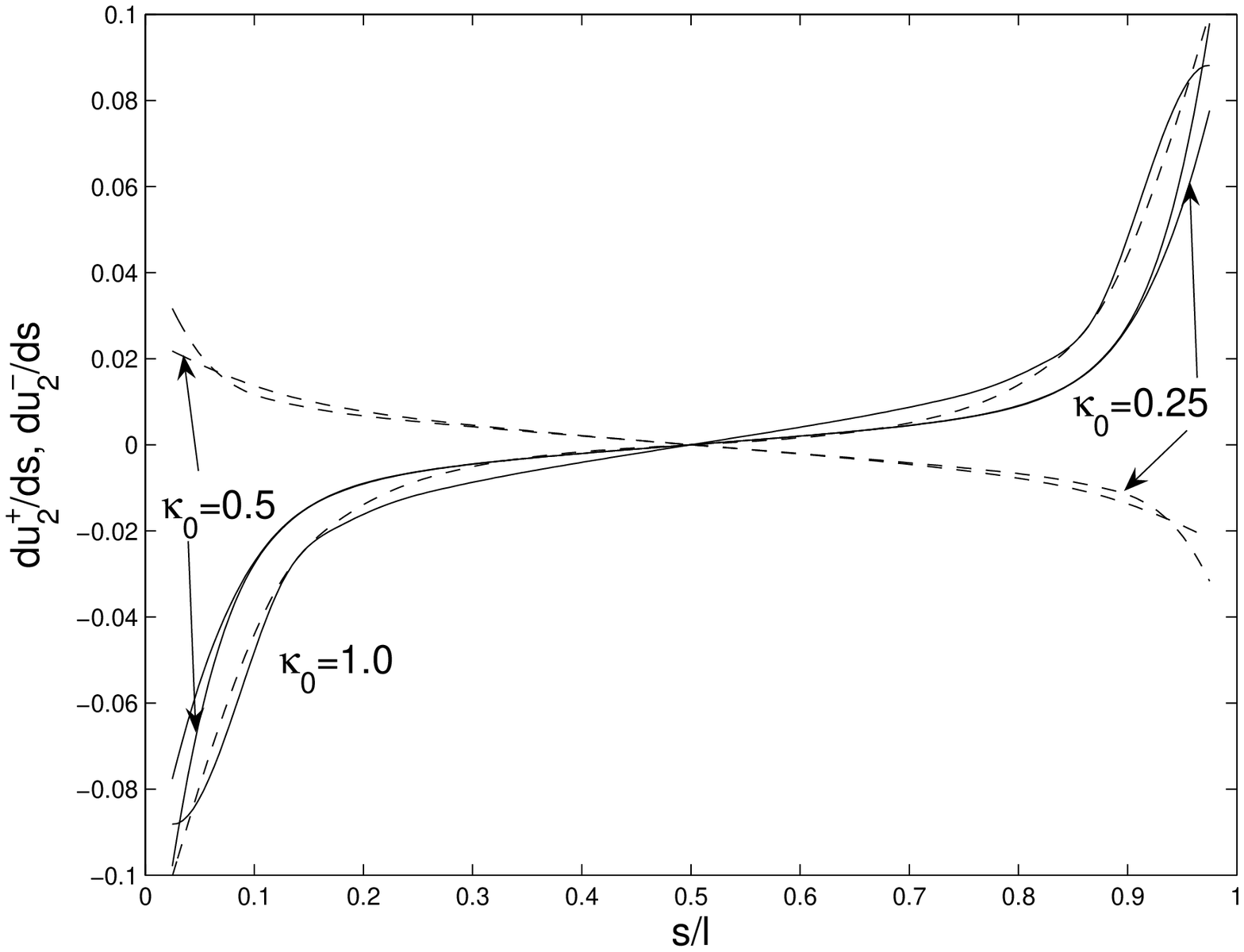}}
	\caption{Graphs of stresses $\sigma^{\pm}_n$, $\tau^{\pm}_n$ and derivatives of displacements $du_1^{\pm}/ds$, $du_2^{\pm}/ds$ for $\sigma_1^{\infty}=1$, $\sigma_2^{\infty}=1$.}
	\label{fig6}
\end{figure}

It has been determined above that the derivative of displacement $du_1^{\pm}/ds$ and the shear stress $\tau_n^{\pm}$ may possess logarithmic singularities at the crack tips. Assume that near the tip of the crack corresponding to the value $s=0$ of the parameter $s$ the functions $du_1^{\pm}/ds$ and $\tau_n^{\pm}$ have the following asymptotic expansions:
$$
\frac{du_1}{ds}=A_1\ln s+O(1),\,\,\, \tau_n=A_2\ln s+O(1).
$$
The graphs of the dependence of the singularity coefficients $A_1$ and $A_2$ on the proportionality constant $\gamma_1$ are presented on the fig. \ref{fig4} for a horizontal ($\sigma_1^{\infty}=1$, $\sigma_2^{\infty}=0$) and a vertical ($\sigma_1^{\infty}=0$, $\sigma_2^{\infty}=1$) stretching at infinity. The numerical results corroborate that the derivative $du_2^{\pm}/ds$ and the stress $\sigma_n^{\pm}$ are bounded at the crack tips.

\begin{figure}[ht]
	\centering
		\scalebox{0.35}{\includegraphics{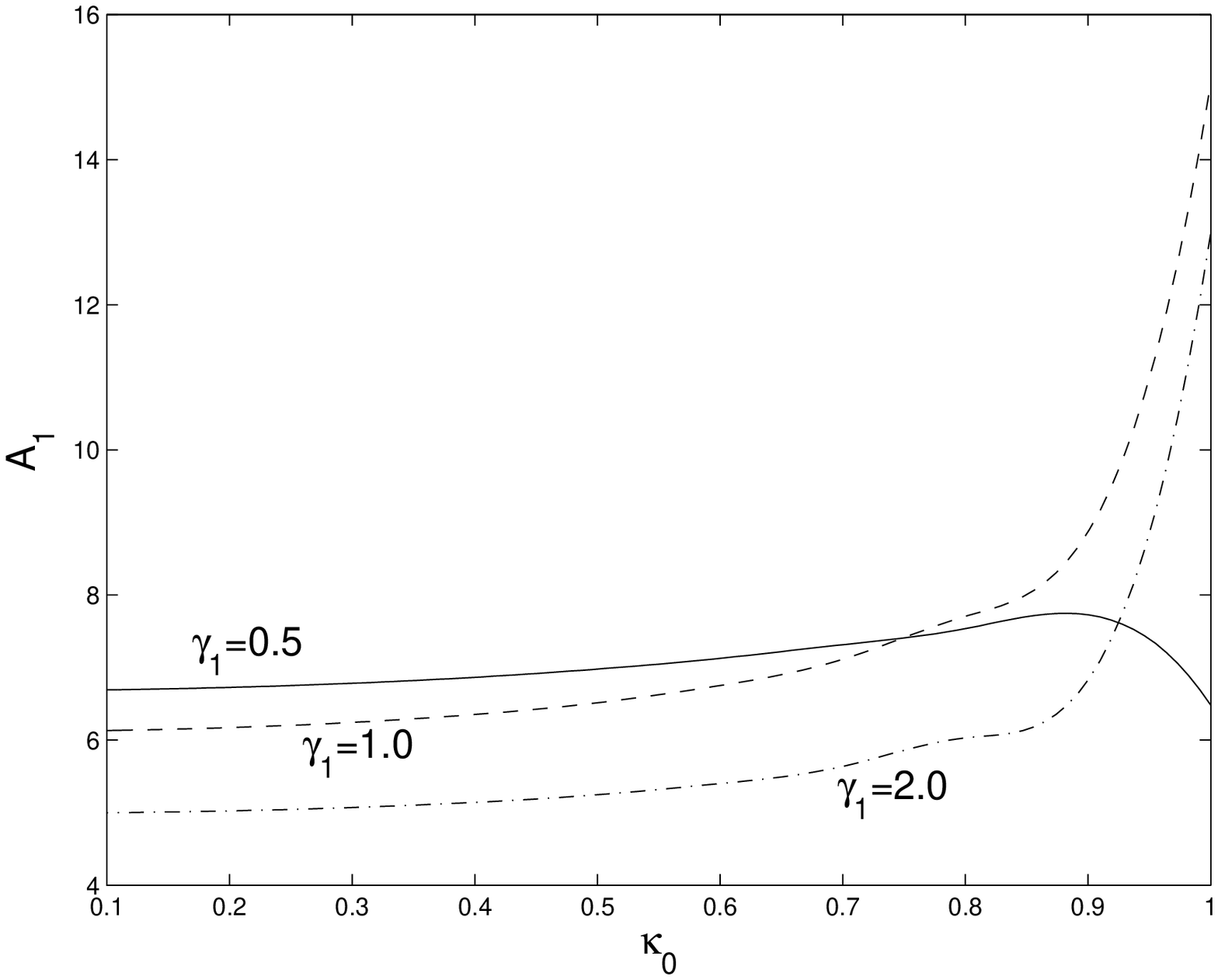}}
		\scalebox{0.35}{\includegraphics{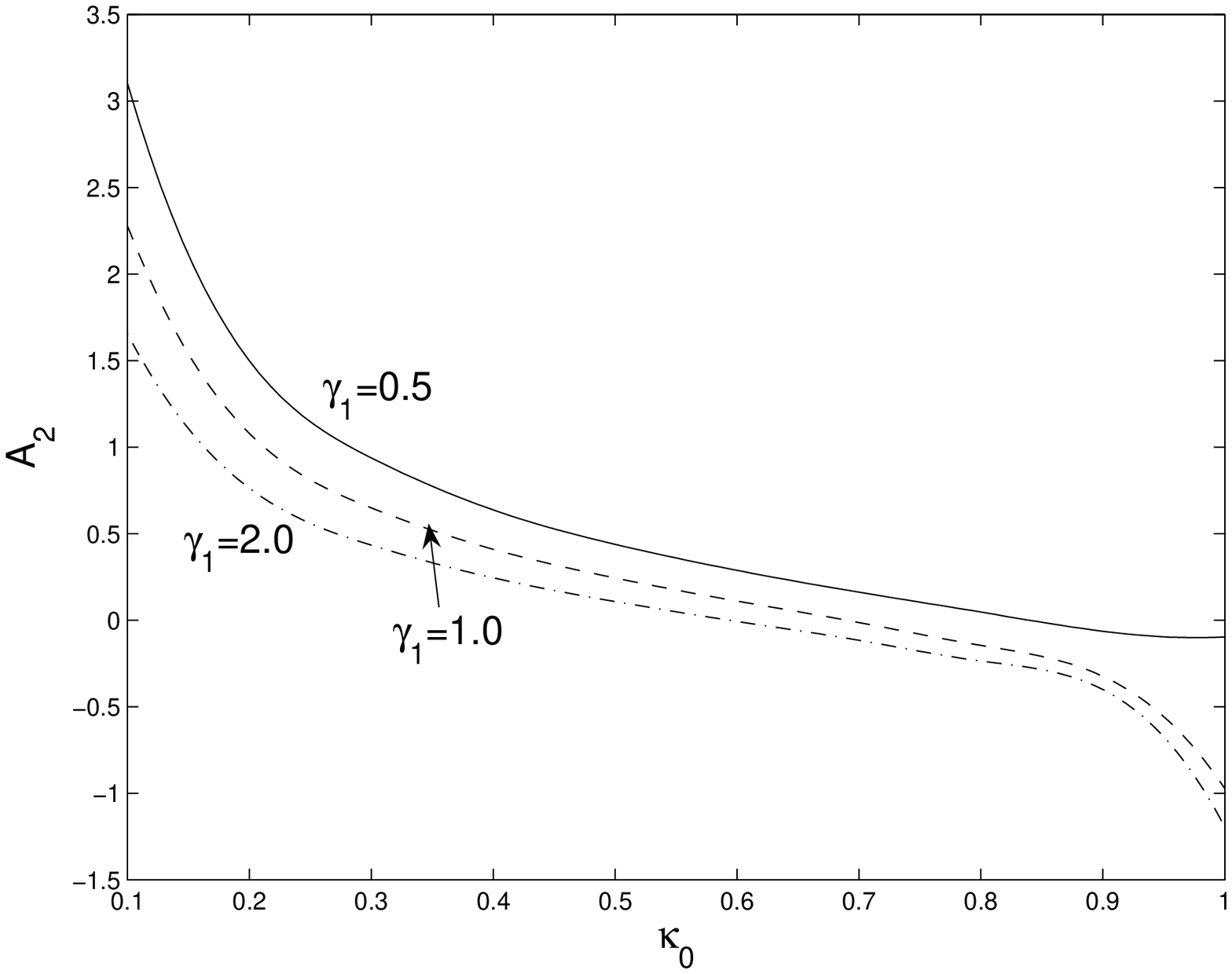}}
			\caption{Graphs of the dependence of the singularity coefficients $A_1$ and $A_2$ on the curvature $\kappa_0$.}
	\label{fig7}
\end{figure}

The graphs of the dependence of the maximal and minimal crack openings on the parameter $\gamma_1$ are presented on the fig. \ref{fig5}. It can be seen that the minimal opening can be negative which suggests that the edges of the crack are pressed together shut in some areas of the crack. To model such behavior, one must allow for contact of the crack surfaces in the boundary conditions with the ends of the contact set being new parameters that must be found as part of the solution procedure.

It can be seen that the dependence of the quantities in the figs. \ref{fig4} and \ref{fig5} on the parameter $\gamma_1$ is highly nonlinear. The consideration of the graphs in the fig. \ref{fig4} and fig. \ref{fig5} shows that the extremal values of all those quantities occur for the same value of $\gamma_1=\gamma_1^{ext}$. The physical meaning and significance of the value $\gamma_1^{ext}$ is immediately unclear and presents an interesting question to study. Since this value appears to be different for different loadings at infinity it can be assumed that $\gamma_1^{ext}$ is dependent on the loading at infinity and not only on the mechanical properties of the plate material and crack geometry. 

The graphs of the dependence of the stresses $\sigma_n^{\pm}$, $\tau_n^{\pm}$ and the derivatives of the displacements $du_1^{\pm}/ds$, $du_2^{\pm}/ds$ on the parameter $s$ are presented on the fig. \ref{fig6} for the circular arc located in the upper half plane between the points $z=\pm 1$ for the three values of the arc curvature: $\kappa_0=0.25$, $\kappa_0=0.5$ and $\kappa_0=1.0$ (fig. \ref{fig5a}). The uniform stretching ($\sigma_1^{\infty}=1$, $\sigma_2^{\infty}=1$) is applied at infinity. The values of other parameters are $\mu=60$ and $\gamma_1=1.0$. As before, the dashed line corresponds to the upper side of the crack (denoted by ``$+$") and the solid line corresponds to the lower side of the crack (denoted by ``$-$").

The graphs of the dependence of the coefficients $A_1$ and $A_2$ on the curvature of the circular arc $\kappa_0$ are presented on the fig. \ref{fig7} for the following values of the parameters $\sigma_1^{\infty}=1$, $\sigma_2^{\infty}=1$ and $\mu=60$. Three different values $\gamma_1=0.5$, $\gamma_1=1.0$ and $\gamma_1=2.0$ of the proportionality coefficient $\gamma_1$ have been considered.

\section{Conclusions}

In this paper a new approach to the fracture modeling first introduced in \cite{OhWaltonSlattery2006}, \cite{SendovaWalton2010} has been extended to the smooth curvilinear cracks of arbitrary shape. The main idea of the approach is to correct the bulk material behavior in the neighborhood of the crack to account for effects of long-range forces from the adjacent phases. Namely, it is assumed that the surface of the crack is endowed with the excess properties dependent on the mean curvature of the crack. 

Due to the fact that the behavior of the material in the bulk is subject to the equations of classical linear elasticity, it is possible to apply the wide range of techniques of LEFM to solve the problem. In this paper the following approach has been taken. The stresses and derivatives of the displacements are expressed through two analytic outside of the crack functions $\Phi(z)$, $\Psi(z)$ (complex potentials) with the help of Muskhelishvili formulas \cite{Mus1963}. The complex potentials can be further written in terms of two functions $g'(s)$ and $q(s)$ by using the integral representations of Savruk \cite{Savruk1981}. Combining these techniques with the new boundary condition which includes the curvature-dependent surface tension, leads to the system of two Cauchy-type singular integro-differential equations with respect to the functions $g'(s)$ and $q(s)$. This system can be further reduced to the system of two weakly singular Fredholm integral equations of the second type and, hence, possesses a unique solution except for countably many values of the parameters.

The numerical results presented in the paper are based on approximation of the unknown function $g'(s)$ by the Taylor polynomials. This leads to the solution of the relatively small systems of linear algebraic equations. It is observed on the examples that the proposed numerical approach provides a good convergence.

It can be seen that there is a principal difference between the case of a straight crack studied in \cite{SendovaWalton2010} and the case of a curvilinear crack considered here. While the stresses and the derivatives of the displacements do not have a classical integrable power singularity of the order $1/2$ as in LEFM, some of them may still have a weaker logarithmic singularity. Thus, more general surface tension models are necessary.


\vspace{.1in}
 
 
\begin{thebibliography}{99}

\bibitem{Abraham2001} {\sc F. F. Abraham.} {\em The atomic dynamic of fracture.} J. Mech. Phys. Solids, 49(2001), pp. 2095-2111.

\bibitem{AbrahamEtal1977} {\sc F. F. Abraham, D. Brodbeck, W. E. Rudge, and X. Xu.} {\em  A molecular dynamics investigation
of rapid fracture mechanics.} J. Mech. Phys. Solids, 45(1997), pp. 1595-1619.

\bibitem{AbrahamGao2000} {\sc F. F. Abraham and H. Gao.} {\em How fast can cracks propagate?} Phys. Rev. Lett., 84(2000), pp. 3113-
3116.

\bibitem{England1971} {\sc A.H. England,} {\em Complex Variable Methods in Elasticity,} Wiley, 1971.

\bibitem{Eshelby1956} {\sc J.D. Eshelby.} {\em The continuum theory of lattice defects.} Progress in Solid State Physics, 3(1956), pp. 79-144.

\bibitem{FinebergGross1991} {\sc J. Fineberg, S. P. Gross, M. Marder, and H. L. Swinney.} {\em Instability in dynamic fracture.}
Phys. Rev. Lett., 67(1991), pp. 457-460.

\bibitem{FometheMaugin1998} {\sc A. Fomethe and G. A. Maugin.} {\em On the crack mechanics of hard ferromagnets.} Internat. J.
Non-Linear Mech., 33(1998), pp. 85-95.

\bibitem{Gakhov1990} {\sc F.D. Gakhov}, {\em Boundary value problems,} Dover Publications, New York, 1990.

\bibitem{Gurtin1995} {\sc M.E. Gurtin.} {\em The nature of configurational forces.} Arch. Rational Mech. Anal., 131(1995), pp. 67-100. 

\bibitem{GurtinPodio1996} {\sc M.E. Gurtin and P. Podio-Guidugli.} {\em Configurational forces and the basic laws for crack
propagation.} J. Mech. Phys. Solids, 44(1996), pp. 905-927.

\bibitem{GurtinPodio1998} {\sc M.E. Gurtin and P. Podio-Guidugli.} {\em Configurational forces and a constitutive theory for crack
propagation that allows for kinking and curving.} J. Mech. Phys. Solids, 46(1998), pp. 1343-1378.

\bibitem{GurtinShvartsman1997} {\sc M. E. Gurtin and M. M. Shvartsman.} {\em Configurational forces and the dynamics of planar cracks
in three-dimensional bodies.} J. Elasticity, 48(1997), pp. 167-191.

\bibitem{HaBobaru2010a} {\sc Y. D. Ha and R. Bobaru.} {\em Studies of dynamic crack propagation and crack branching with
peridynamics.} Int. J. Frac., 162(2010), pp. 229-244.

\bibitem{HaBobaru2011} {\sc Y. D. Ha and R. Bobaru.} {\em Characteristics of dynamic brittle fracture captured with peridynamics.}
Engr. Frac. Mech., 78(2011), pp. 1156-1168.

\bibitem{HollandMarder1997} {\sc D. Holland and M. P. Marder.} {\em Ideal brittle fracture of silicon studied with molecular dynamics.}
Phys. Rev. Lett., 80(1997), pp. 746-749.

\bibitem{MarderGross1995} {\sc M. Marder and S. Gross.} {\em Origin of crack tip instabilities.} J. Mech. Phys. Solids, 43(1995), pp. 1-48.

\bibitem{MauginTrimarco1995} {\sc G. A. Maugin and C. Trimarco.} {\em Dissipation of configurational forces in defective elastic solids.}
Arch. Mech. (Arch. Mech. Stos.), 47(1995), pp. 81-99.

\bibitem{MikhPros1986} {\sc S. G. Mikhlin, S. Pr\"ossdorf}, {\em Singular integral operators,} Springer, Berlin, 1986. (Translated from the German by Albrecht Böttcher and Reinhard Lehmann.)

\bibitem{Mus1963} {\sc N. I. Muskhelishvili,} {\em Some basic problems of the mathematical theory of elasticity; fundamental equations, plane theory of elasticity, torsion, and bending,} Noordhoff International Publishing, Groningen, 1963.

\bibitem{OhWaltonSlattery2006} {\sc E.-S. Oh, J. R. Walton, and J. C. Slattery.} {\em A Theory of Fracture Based Upon an Extension of Continuum Mechanics to Nanoscale.} J. Appl. Mech., 73(2006), pp. 792-798.

\bibitem{Savruk1981} {\sc M.P. Savruk,} {\em Two-dimensional problems of elasticity for cracked solids,} Naukova dumka, Kiev (in Russian), 1981.

\bibitem{SendovaWalton2010} {\sc T. Sendova,  J.R. Walton.} {\em A new approach to the modeling and analysis of fracture through extension of continuum mechanics to the nanoscale.} Mathematics and Mechanics of Solids, 15(2010), pp. 368-413. 

\bibitem{Silling2000} {\sc S. A. Silling.} {\em Reformulation of elasticity theory for discontinuities and long-range forces.} J.
Mech. Phys. Solids, 48(2000), pp. 175-209.

\bibitem{Slattery1990} {\sc J.C. Slattery,} {\em Interfacial Transport Phenomena,} Springer-Verlag, New York, 1990.

\bibitem{SlatteryEtAl2004} {\sc J. C. Slattery, E. S. Oh, and Kai-Bin Fu.} {\em Extension of continuum mechanics to the nanoscale.}
Chem. Eng. Sci., 59(2004), pp. 4621-4635.

\bibitem{SlepyanEtAl1999} {\sc L. I. Slepyan, M. V. Ayzenberg-Stepanenko, and J. P. Dempsey.} {\em A lattice model for viscoelastic
fracture.} Mech. Time-Dep. Mater., 3(1999), pp. 159-203.

\bibitem{Sokolnikoff1956} {\sc I.S. Sokolnikoff,} {\em Mathematical Theory of Elasticity,} McGraw-Hill, 1956. 

\bibitem{Spector2002} {\sc J. Sivaloganathan and S. J. Spector.} {\em On cavitation, configurational forces and implications
for fracture in a nonlinearly elastic material.} J. Elasticity, 67(2002), pp. 25-49.

\bibitem{SwadenerEtAl2002} {\sc J. G. Swadener, M. I. Baskes, and M. Nastasi.} {\em Molecular dynamics simulation of brittle fracture
in silicon.} Phys. Rev. Lett., 89(2002), 085503.

\bibitem{TadmoretAl1996} {\sc E. B. Tadmor, R. Philips, and M. Ortiz.} {\em Mixed atomistic and continuum models of deformation
in solids.} Langmuir, 12(1996), pp. 4529-4534.

\bibitem{XiaoBelytschko2004} {\sc S. P. Xiao and T. Belytschko.} {\em A bridging domain method for coupling continua with molecular
dynamics.} Comp. Meth. Appl. Mech. \& Engr., 193(2004), pp. 1645-1669.


\end{thebibliography}
\end{document}